\documentclass[a4paper, 11pt]{amsart}
\usepackage{amssymb}
\usepackage{amsthm}
\usepackage{fullpage}
\usepackage{graphicx}
\usepackage{subfig}
\usepackage{amsmath,amscd}
\usepackage{xypic}
\usepackage{color}
\usepackage{float}
\usepackage{vmargin}
\usepackage[colorlinks,citecolor=blue,filecolor=black,linkcolor=blue,urlcolor=black]{hyperref}
\usepackage[all,cmtip]{xy}

\newcommand{\Z}{\ensuremath{\mathbb Z}}

\newcommand{\R}{\ensuremath{\mathbb R}}

\theoremstyle{plain}
\newtheorem{thm}{Theorem}[section]
\newtheorem*{thm*}{Theorem}

\newtheorem*{cor*}{Corollary}
\newtheorem*{prop*}{Proposition}

\newtheorem*{lemma*}{Lemma}
\newtheorem*{claim*}{Claim}

\theoremstyle{definition}

\newtheorem*{exmp*}{Example}
\newtheorem*{defn*}{Definition}
\newtheorem*{rem*}{Remark}
\newtheorem*{note*}{Note}

\begin{document}

\title{Kirby diagrams and the Ratcliffe-Tschantz hyperbolic 4-manifolds} 

\author{Hemanth Saratchandran}
\email{hemanth.saratchandran@maths.ox.ac.uk}
\date{\today}

\maketitle 
\parskip=0.2cm
\parindent=0.0cm

\begin{abstract}
We show how to construct a Kirby diagram for a large class of finite volume hyperbolic 4-manifolds constructed by J. Ratcliffe and S. Tschantz.
\end{abstract}

\tableofcontents

\section{Introduction}

The primary aim of this paper is to initiate a study of finite volume hyperbolic 4-manifolds using techniques from 4-manifold theory. Although there have
been many constructions of finite volume hyperbolic 4-manifolds (see \cite{kolpakov}, \cite{ratcliffe}, \cite{slavich}) their exploration using invariants coming from 4-manifold theory remains
a largely unexplored topic. One of the most powerful and general techniques for studying various existence questions about invariants associated to a
4-manifold comes from the theory of Kirby diagrams and the associated Kirby calculus. This approach to probing the internal structure of a 4-manifold
gained prominence in the 70's and was turned in to a fine art soon after by many topologists. The focus in this paper is to explicitly construct
Kirby diagrams for a certain collection of finite volume hyperbolic 4-manifolds. One of the largest census of finite volume hyperbolic 4-manifolds was
constructed by J. Ratcliffe and S. Tschantz in their paper \cite{ratcliffe}. They construct a total of 1171 distinct isometry classes of finite volume hyperbolic
4-manifolds, each with Euler characteristic one. The fundamental domain for the assortment of manifolds they construct is the 24-cell, a four dimensional
self dual polyhedron. Using some simple symmetry, and the self dual nature of this polyhedron we show how to construct a Kirby diagram for each one
of the 4-manifolds constructed by Ratcliffe and Tschantz. We explain how to identify the one, two and three handles of each of the 4-manifolds giving various
pictures along the way. In a future paper we will show how to use the constructions in this paper to produce explicit examples of four dimensional
hyperbolic link complements.



\section*{Acknowledgements}
The author wishes to acknowledge Andras Juhasz, Panos Papazoglou, John Parker and Marc Lackenby for comments and corrections on an earlier
draft of this work.

\section{Dualising the 24-cell}\label{dual_24}

In this section we will explain to the reader how the 24-cell is constructed, and how we can visualise the boundary of the 24-cell in $\R^3$. 
We will not be going in to details about the structure of the 24-cell as our primary focus will be on handle decompositions.
For more background information on the 24-cell we refer the reader to \cite{cox} chap.4. The paper \cite{kerckhoff} has some background on the
24-cell with some nice pictures depicting various edges and faces of the 24-cell.   

Let $S_{(*,*,*,*)}$ denote a sphere of radius 1 centred at a point in $\R^4$ whose coordinates have two $\pm 1$'s and whose other two coordinates are both zero.
For example $S_{(+1,+1,0,0)}$ denotes the sphere of radius 1 centred at the point $(1,1,0,0)$ in $\R^4$. If we let $\mathbb{H}^4$ denote the ball model 
of hyperbolic 4-space, then we find that all the spheres $S_{(*,*,*,*)}$ intersect the sphere at infinity orthogonally. 
This implies that each such sphere determines a hyperplane in $\mathbb{H}^4$. 
If we let $Q_{(*,*,*,*)}$ denote the corresponding half-space that contains the origin, and then take the intersection of all such half-spaces, we find that we obtain
a 24-sided polyhedron $P$ in $\mathbb{H}^4$. This polyhedron is known as the hyperbolic 24-cell and we will denote it by $P$. It is a four dimensional self dual
polyhedron.
We will denote the side of $P$ that lies on the sphere
$S_{(*,*,*,*)}$ also by $S_{(*,*,*,*)}$. All the dihedral angles of $P$ are $\pi/2$ and it has 24 vertices which are all ideal vertices. We can explicitly describe each ideal
vertex: We have 8 vertices of the form $v_{(\pm 1, 0,0,0)} = (\pm 1,0,0,0)$, $v_{(0, \pm 1, 0,0)} = (0,\pm 1,0,0)$, $v_{(0,0, \pm 1, 0)} = (0,0,\pm 1,0)$, 
$v_{(0,0,0, \pm)} = (0,0,0,\pm 1)$ plus 16 vertices of the form 
$v_{(\pm 1/2, \pm 1/2, \pm 1/2, \pm 1/2)} = (\pm 1/2, \pm 1/2, \pm 1/2, \pm 1/2)$. Finally let us mention that it has twenty four codimension 1 sides, 
ninety six codimension 2 sides, and ninety six codimension 3 sides.

In our study of the 24-cell it will be very useful for us to know which sides intersect in codimension 2 faces and which vertices lie on a particular side.
Two sides of $P$ intersect in $\mathbb{H}^4$ if their identifying spheres have coordinates that have equal nonzero entries in one place and the remaining nonzero entries lie in different positions. For example the sides $S_{(+1,0,-1,0)}$ and $S_{(+1,+1,0,0)}$ intersect, however the sides $S_{(+1,+1,0,0)}$ and $S_{(0,0,-1,-1)}$ do not. In the case of vertices 
with one nonzero $\pm 1$ (those from the first 8 described above) we find that such a vertex lies on a side if the side has an equal nonzero entry in the same
position as the vertex. For example the vertex $v_{(-1,0,0,0)}$ lies on the side $S_{(-1,0,0,+1)}$ but not on the side $S_{(0,+1,+1,0)}$. In the case of a vertex in the group
$v_{(\pm 1/2,\pm 1/2,\pm 1/2,\pm 1/2)}$ we find that such a vertex lies on a side provided the sign of the nonzero entries of the side coincides with the sign of the entries of the vertex in the same position.
For example the vertex $v_{(+1/2, +1/2, -1/2, +1/2)}$ lies on the side $S_{(+1,0,-1,0)}$ but does not lie on the side $S_{(-1,+1,0,0)}$. We also mention that
when two sides of $P$ intersect they do so at right angles. 

The hyperbolic manifolds that we will be considering will be obtained by a side pairing identification of the 24-cell. In this regard it will be very helpful
if we had some way to visualise the 24-cell. The 24-cell is a 4-dimensional self dual polyhedron, hence its boundary is 3-dimensional. We will now outline a procedure that
will allow us to visualise this boundary. This will prove useful when we start looking at a Kirby diagram of a particular manifold obtained from a side pairing
transformation of the 24-cell.

Each side of the polyhedron $P$ corresponds to a sphere of the form $S_{(*,*,*,*)}$, which itself is identified by the co-ordinate $(*,*,*,*)$. We can then radially project this point to the boundary $\partial{B^4} = S^3$. For example, if we are focusing on the side corresponding
to the sphere $S_{(+1,0,+1,0)}$, then the radial projection of the point $(1,0,1,0)$ is $(1/\sqrt{2}, 0, 1/\sqrt{2}, 0)$. Thus each side in $P$ uniquely corresponds to
a point in $S^3$.

Using the Mobius transformation that transforms the ball $B^4$ onto the upper half-plane $\R^4_{x_4 > 0}$ we obtain a map from $S^3$ onto $\R^3 \cup \{\infty\}$. We will
use this map to transfer information about the boundary of $P$ into $\R^3 \cup \{\infty\}$. 

We define the map:
\[ \phi : S^3 \rightarrow \R^4 \cup \{\infty\} \]
by 
\[ \phi(x_1, x_2, x_3, x_4) = (0, 0, 0, 1) + \frac{2}{x_1^2 + x_2^2 + x_3^2 + (x_4 - 1)^2}(x_1, x_2, x_3, x_4 - 1). \]
This is just the usual Mobius transformation from $\overline{B^4}$ to  $\R^4_{x_4 \geq 0}\cup \infty$ restricted to the boundary sphere $S^3$. 
Using the above map we can map each point in $S^3$ that corresponds to a sphere onto a point in $\R^3$. For example if we take the sphere $S_{(+1,0,+1,0)}$, then we know that
the associated point on $S^3$ has coordinates $(1/\sqrt{2}, 0, 1/\sqrt{2}, 0)$. Applying $\phi$ to this point we get the point  $(1/\sqrt{2}, 0, 1/\sqrt{2})$. Doing this
for all the points corresponding to the sides of $P$ we get a collection of points in $\R^3$.

Suppose we have two sides of the polyhedron $P$ call them $S_1$ and $S_2$ that intersect along a codimension two face (a codimension two face is sometimes called a ridge). Associated to these sides we have points $x_1$ and $x_2$ respectively that reside on $S^3$. We can then choose a path in $S^3$ from $x_1$ to $x_2$ that when projected to $P$ will intersect the codimension 2 side 
$S_1 \cap S_2$ transversely in one point. Let us note that there will be many paths we can choose between $x_1$ to $x_2$, but if we give $S^3$ the round metric induced from the
standard Euclidean inner product on $\R^4$, then there will be a unique length minimising geodesic segment joining $x_1$ to $x_2$. It is this segment that we choose as our path, and from here on
in when we speak of such a path between two such points we keep it fixed that it is the geodesic segment we are taking as the path.  
We can think of this path in $S^3$ as representing the intersecting face $S_1 \cap S_2$. For example if we took the sides
$S_{(+1,+1,0,0)}$ and $S_{(+1,0,+1,0)}$ then we know that the corresponding points in $S^3$ are given by the points $(1/\sqrt{2}, 1/\sqrt{2}, 0, 0)$ and 
$(1/\sqrt{2}, 0, 1/\sqrt{2}, 0)$. The corresponding path is taken to be the geodesic segment between these two points. If we then use our map $\phi$ these paths
between such points will give us a path in $\R^3$ between the images of these points. In this way we view the image path in $\R^3$ as representing a codimension 2 face in
the polyhedron $P$, its endpoints correspond to the sides in $P$ that intersect to give this codimension 2 face. We mention that since the path we are taking
in $S^3$ that corresponds to the intersection $S_1 \cap S_2$ is a geodesic segment (with respect to the round metric on $S^3$) the image will in general not
be a straight line. However, the image path will be homotopic (fixing endpoints) to a straight line and when drawing pictures we will always take the image
paths to be straight line segments.

Thus far we have explained how to visualise the faces and edges (codimension 2 faces) of the polyhedron $P$, before we draw these edges let us explain how we can
also visualise the codimension 3 faces. Fix three sides $S_1$, $S_2$ and $S_3$, and suppose they all triply intersect in a codimension 3 face. 
Let $x_1$, $x_2$, $x_3$ denote the respective points in $S^3$. As we explained above each pair of intersections defines a geodesic segment in $S^3$ whose endpoints
are precisely the points corresponding to the intersecting side. Thus associated to the triply intersecting sides $S_1$, $S_2$ and $S_3$ we have a geodesic
triangle in $S^3$ whose vertices are $x_1$, $x_2$ and $x_3$, and whose sides correspond to the pair of intersections of $S_1$, $S_2$ and $S_3$. This geodesic
triangle in $S^3$ will get mapped to a triangle in $\R^3$ via $\phi$. In general this triangle will not be made up of straight line segments, however when
we draw such triangles we will always take straight line segments. This tells us that codimension 3 faces correspond to triangles made up of
paths that correspond to codimension 2 faces.

Putting together the above information gives us a way to picture the boundary of $P$ in $\R^3$. What we are really visualising is the dual polyhedron of $\partial{P}$, but 
since $P$ is self dual this amounts to visualising the boundary of $P$.

\centerline {\includegraphics[width=7cm, height=7cm]{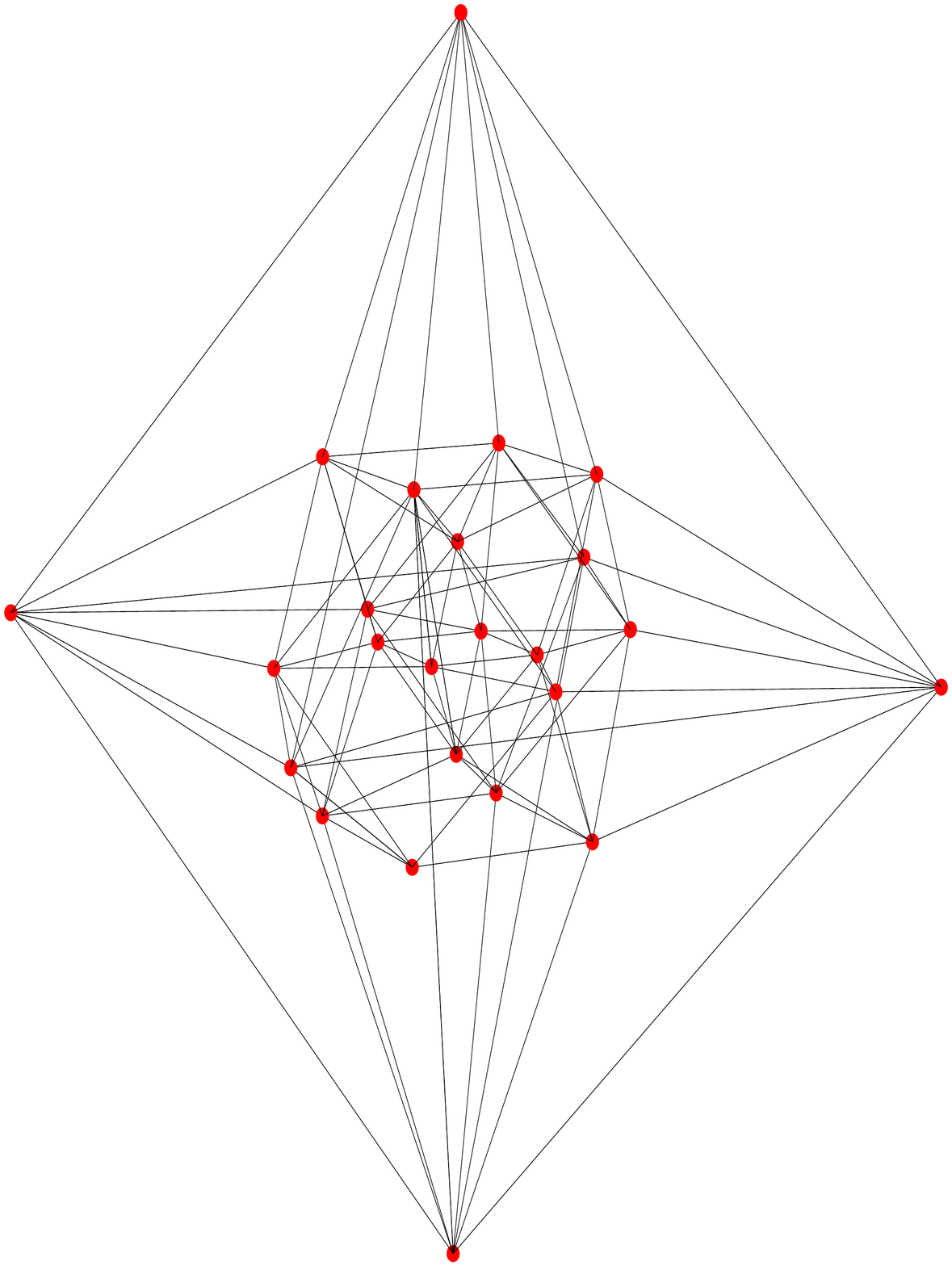}}

The red vertices in the above diagram correspond to the codimension 1 sides of $P$ and recall that we obtained these by picking the centre point of each sphere, radially
projecting it to $S^3$, then mapping it to $\R^3$ using the Mobius transformation $\phi$. The following table shows the co-ordinates of each of the red vertices and their
corresponding sphere. \\

\begin{tabular}{|l|l||l|l|}
 $S_{(+1,+1,0,0)}$ & $(\frac{1}{\sqrt{2}}, \frac{1}{\sqrt{2}}, 0)$  &  $S_{(-1,+1,0,0)}$ & $(\frac{-1}{\sqrt{2}}, \frac{1}{\sqrt{2}}, 0)$  \\
 $S_{(+1,-1,0,0)}$ & $(\frac{1}{\sqrt{2}}, \frac{-1}{\sqrt{2}}, 0)$ &  $S_{(-1,-1,0,0)}$  & $(\frac{-1}{\sqrt{2}}, \frac{-1}{\sqrt{2}}, 0)$ \\
 $S_{(+1,0,+1,0)}$ &  $(\frac{1}{\sqrt{2}}, 0, \frac{1}{\sqrt{2}})$ &  $S_{(+1,0,-1,0)}$ & $(\frac{1}{\sqrt{2}}, 0, \frac{-1}{\sqrt{2}})$ \\
 $S_{(-1,0,+1,0)}$ &  $(\frac{-1}{\sqrt{2}}, 0, \frac{1}{\sqrt{2}})$ & $S_{(-1,0,-1,0)}$ & $(\frac{-1}{\sqrt{2}}, 0, \frac{-1}{\sqrt{2}})$ \\ 
 $S_{(0,+1,+1,0)}$ &  $(0, \frac{1}{\sqrt{2}} ,\frac{1}{\sqrt{2}})$ &  $S_{(0,-1,-1,0)}$ &  $(0, \frac{-1}{\sqrt{2}} ,\frac{-1}{\sqrt{2}})$ \\ 
 $S_{(0,+1,-1,0)}$ &   $(0, \frac{1}{\sqrt{2}} ,\frac{-1}{\sqrt{2}})$ &  $S_{(0,-1,+1,0)}$ &  $(0, \frac{-1}{\sqrt{2}} ,\frac{1}{\sqrt{2}})$ \\ 
 $S_{(+1,0,0,+1)}$ &  $(1 + \sqrt{2}, 0, 0)$ &  $S_{(-1,0,0,-1)}$ & $(1 - \sqrt{2}, 0, 0)$ \\ 
 $S_{(+1,0,0,-1)}$ &  $(-1 + \sqrt{2}, 0, 0)$ &  $S_{(-1,0,0,+1)}$ & $(-1 - \sqrt{2}, 0, 0)$ \\
 $S_{(0,+1,0,+1)}$ &  $(0, 1 + \sqrt{2}, 0)$ &  $S_{(0,-1,0,+1)}$ &  $(0, -1 - \sqrt{2}, 0)$ \\
 $S_{(0,+1,0,-1)}$ &   $(0, -1 + \sqrt{2}, 0)$ &  $S_{(0,-1,0,-1)}$ &  $(0, 1 - \sqrt{2}, 0)$  \\
 $S_{(0,0,+1,+1)}$ &   $(0, 0, 1 + \sqrt{2})$ &  $S_{(0,0,+1,-1)}$ &  $(0, 0, -1 + \sqrt{2})$   \\
 $S_{(0,0,-1,+1)}$ &   $(0, 0, -1 - \sqrt{2})$  &  $S_{(0,0,-1,-1)}$ &  $(0, 0, 1 - \sqrt{2})$  
\end{tabular}

\section{Preliminaries on the Ratcliffe-Tschantz hyperbolic 4-manifolds}

In their paper \cite{ratcliffe} Ratcliffe and Tschantz construct several examples of non-compact hyperbolic 4-manifolds of finite volume. We briefly describe how this
construction process works. The interested reader is referred to page 109 of their paper \cite{ratcliffe} for more details.

Recall, a real $(n+1) \times (n+1)$ matrix $A$ is said to be \textit{Lorentzian} if $A$ preserves the Lorentzian inner product:
\[ x \cdot y := x_1y_1 + \ldots + x_ny_n  - x_{n+1}y_{n+1}. \]
The hyperboloid model of hyperbolic n-space is the metric space 
\[ \mathbb{H}^n := \{ x \in \R^{n+1} : x \cdot x = -1 \text{ and } x_{n+1} > 0 \} \]
with the metric $d$ defined by the equation:
\[ cosh(d(x,y)) := -x \cdot y .\] 
A Lorentzian $(n+1) \times (n+1)$ matrix $A$ is said to be positive if $A$ maps $\mathbb{H}^n$ to $\mathbb{H}^n$. The isometries of $\mathbb{H}^n$ are the positive
Lorentzian matrices. \\
Let $\Gamma^n$  be the group of $(n+1) \times (n+1)$ positive Lorentzian matrices with integer entries. The group $\Gamma^n$ is an infinite discrete subgroup of the group $O(n, 1)$ of
Lorentzian $(n+1) \times (n+1)$ matrices. The principal congruence two subgroup of $\Gamma^n$ is the group $\Gamma^n_2$ of all matrices in $\Gamma^n$ that are congruent
to the identity matrix modulo two.

In their paper \cite{ratcliffe} Ratcliffe and Tschantz classify all the hyperbolic space forms $\mathbb{H}^4/\Gamma$ where $\Gamma$ is a torsion free subgroup
of minimal index in the group $\Gamma^4_2$ (they actually do the case $n = 2, 3$ as well but as we are only interested in the $n = 4$ case we will not worry about these
other two). The main theorem takes the form:

\begin{thm}
There are, up to isometry, exactly 1171 hyperbolic space-forms $\mathbb{H}^4/\Gamma$ where $\Gamma$ is a torsion free subgroup of minimal index in the group
$\Gamma^4_2$. Only 22 of these manifolds are orientable.
\end{thm}

The proof of this theorem which leads to the construction of the 1171 hyperbolic space-forms $\mathbb{H}^4/\Gamma$ involves finding suitable side pairings of the 24-cell that will give rise to hyperbolic 4-manifolds. The interested reader can consult p.111 of \cite{ratcliffe} for the details of the proof.
The actual side-pairings they use will be important for our understanding of a Kirby diagram, so we go through in some detail exactly how they code their
side pairings. 

All their side pairings are of the form $rk$, where $k$ is a diagonal matrix with diagonal entry taking the form $(\pm 1, \pm 1, \pm 1, \pm 1)$, which is to be interpreted as 
a composition of reflections in the coordinate planes $x_i = 0$ for $1 \leq i \leq 4$. In order to
identify a particular element $k$ it suffices to simply give the diagonal  $(\pm 1, \pm 1, \pm 1, \pm 1)$, where a $- 1$ in some position tells us that we reflect in the coordinate 
corresponding to that position, and $+ 1$ tell us that we do nothing. For example $(-1, +1, -1, +1)$ is the composition of the reflection in the hyper-planes $x_1 = 0$ followed by 
reflection in $x_3 = 0$. When we want to speak of a particular side pairing $rk$ and make reference to its $k$-part we will always write the $k$-part in the form
$k_{_{(\pm 1, \pm 1, \pm 1, \pm 1)}}$, the $k$ is to remind us that we are only dealing with the $k$ part of $rk$ and the subscript $(\pm 1, \pm 1, \pm 1, \pm 1)$ tells us
what the diagonal is.
Ratcliffe and Tschantz develop a coding system for the $k$-part of each side pairing which we now describe. 

The polyhedron $P$ has 24 three dimensional sides, since a side pairing transformation must identify pairs of sides we need to give twelve transformations. We will
denote these transformations by the letters $a, b, \ldots , k, l$.
We then group the letters $a, b, \ldots ,l$ and the sides of $P$ into the following groups:
\[ \{a, b, S_{(\pm 1,\pm 1, 0,0)}\}, \{c, d, S_{(\pm 1, 0, \pm 1, 0)}\}, \{e, f, S_{(0,\pm 1, \pm 1,0)}\}, \]
\[\{g, h, S_{(\pm 1, 0, 0, \pm 1)}\}, \{i, j, S_{(0,\pm 1, 0, \pm 1)}\}, \{k, l, S_{(0, 0, \pm 1, \pm 1)}\}. \]
In each of the above sets the letters are pairings between the spheres in that set (remember the spheres represent the sides of the 24-cell $P$ by the way we constructed it). We give 
the sides the ordering $(+1, +1) < (+1, -1) < (-1, +1) < (-1, -1)$. So for example taking the first group above our order tells us that 
\[ S_{(+1, +1, 0,0)} <  S_{(+1,-1, 0,0)} < S_{(-1, +1, 0,0)} < S_{(-1, -1, 0,0)}  \]
The first letter always pairs the side with $+1, +1$ to one of the other sides (we will show how to determine this side), and the second letter pairs the next unused side
(with respect to the above ordering) with the last side.
The actual side pairing transformations are encoded by a string of six characters from the set 
\[ \{1,2,3,4,5,6,7,8,9,A,B,C,D,E,F\} \]
one character for each of the above sets. Each character represents a particular $k$-part, and the following table shows the correspondence:

\begin{center}
\begin{tabular}{|l | l | }
\hline	

\textbf{Character} &  \textbf{k-part}  \\ \hline

1 & $k_{_{(-1, +1, +1, +1)}}$  \\ \hline

2 & $k_{_{(+1, -1, +1, +1)}}$  \\ \hline

3 & $k_{_{(-1, -1, +1, +1)}}$  \\ \hline

4 & $k_{_{(+1, +1, -1, +1)}}$ \\ \hline

5 & $k_{_{(-1, +1, -1, +1)}}$ \\ \hline
 
6 & $k_{_{(+1, -1, -1, +1)}}$ \\ \hline 
 
7 & $k_{_{(-1, -1, -1, +1)}}$ \\ \hline 
 
8 & $k_{_{(+1, +1, +1, -1)}}$ \\ \hline  
 
9 & $k_{_{(-1, +1, +1, -1)}}$ \\ \hline   
 
A & $k_{_{(+1, -1, +1, -1)}}$ \\ \hline 
 
B & $k_{_{(-1, -1, +1, -1)}}$ \\ \hline  
 
C & $k_{_{(+1, +1, -1, -1)}}$ \\ \hline   
 
D & $k_{_{(-1, +1, -1, -1)}}$ \\ \hline    
 
E & $k_{_{(+1, -1, -1, -1)}}$ \\ \hline  
 
F & $k_{_{(-1, -1, -1, -1)}}$ \\ \hline

\end{tabular} 
\end{center}

The coding of a particular manifold will take the form of six characters from the above set of characters. For example, a code can look like
\[ \textbf{1428BD} .\]
Each character in the above code tells us what the $k$-part of each pair of transformations in the group
\[ \{(a, b), (c, d), (e, f), (g, h), (i, j), (k, l)\} \]
is. For example for the code \textbf{1428BD} the first character is \textbf{1}, this tells us that for the particular manifold corresponding to this code
the side pairings $a$ and $b$ have $k$-parts given by $k_{_{(-1, +1, +1, +1)}}$. The next character in the code is \textbf{4}, this tells us that the side pairings
$c$ and $d$ have $k$-part  $k_{_{(+1, +1, -1, +1)}}$. Continuing in this way we can determine all the $k$-parts of each of the side pairings $a, b, c, d, \ldots , k, l$.

\[\xymatrixcolsep{5pc}\xymatrix{ S_{(+1,+1,0,0)}  \ar[r]^a_{k_{_{(-1,+1,+1,+1)}}} & S_{(-1,+1,0,0)} } \hspace{2cm} \xymatrix{S_{(+1,-1,0,0)}  \ar[r]^b_{k_{_{(-1,+1,+1,+1)}}} & S_{(-1,-1,0,0)} } \]

\[\xymatrixcolsep{5pc}\xymatrix{ S_{(+1,0,+1,0)}  \ar[r]^c_{k_{_{(+1,+1,-1,+1)}}} & S_{(+1,0,-1,0)} } \hspace{2cm} \xymatrix{S_{(-1,0,+1,0)}  \ar[r]^d_{k_{_{(+1,+1,-1,+1)}}} & S_{(-1,0,-1,0)} } \]

\[\xymatrixcolsep{5pc} \xymatrix{ S_{(0,+1,+1,0)}  \ar[r]^e_{k_{_{(+1,-1,+1,+1)}}} & S_{(0,-1,+1,0)} } \hspace{2cm} \xymatrix{S_{(0,+1,-1,0)}  \ar[r]^f_{k_{_{(+1,-1,+1,+1)}}} & S_{(0,-1,-1,0)} } \]

\[\xymatrixcolsep{5pc} \xymatrix{ S_{(+1,0,0,+1)}  \ar[r]^g_{k_{_{(+1,+1,+1,-1)}}} & S_{(+1,0,0,-1)} } \hspace{2cm} \xymatrix{S_{(-1,0,0,+1)}  \ar[r]^h_{k_{_{(+1,+1,+1,-1)}}} & S_{(-1,0,0,-1)} } \]

\[\xymatrixcolsep{5pc} \xymatrix{ S_{(0,+1,0,+1)}  \ar[r]^i_{k_{_{(-1,-1,+1,-1)}}} & S_{(0,-1,0,-1)} } \hspace{2cm} \xymatrix{S_{(0,+1,0,-1)}  \ar[r]^j_{k_{_{(-1,-1,+1,-1)}}} & S_{(0,-1,0,+1)} } \]

\[\xymatrixcolsep{5pc} \xymatrix{ S_{(0,0,+1,+1)}  \ar[r]^k_{k_{_{(-1,+1,-1,-1)}}} & S_{(0,0,-1,-1)} } \hspace{2cm} \xymatrix{S_{(0,0,+1,-1)}  \ar[r]^l_{k_{_{(-1,+1,-1,-1)}}} & S_{(0,0,-1,+1)} }. \]

The above is to be understood as the letter above the arrow tells you which side of the polyhedron $P$ is being paired to another side, and the $k$-part of that
transformation is on the bottom of the arrow.

We have still not explained what the $r$-part of the side pairing transformations are, they are just reflections in the image side, viewed as a hyperplane in $\mathbb{H}^4$.
For example for the code \textbf{1428BD}, described above, we know that the side pairing transformation $a$ pairs side $S_{(+1,+1,0,0)}$ to side $S_{(-1,+1,0,0)}$, and
has $k$-part given by $k_{_{(-1,+1,+1,+1)}}$. Using the fact that the $r$-part is just reflection in the image side we have that
\[ a := rk_{_{(-1,+1,+1,+1)}} \]
where $r$ is reflection in the side $S_{(-1,+1,0,0)}$.

To summarise, Ratcliffe and Tschantz construct 1171 distinct non-compact finite volume hyperbolic 4-manifolds, 22 of which are orientable and all others non-orientable.
Each such 4-manifold is obtained by giving side pairing transformations for the sides of the 24-cell $P$. The side
pairing transformations for any fixed such manifold are labelled by the letters $a, b, \ldots , k, l$, and each such transformation is given in the form
$rk$, where $r$ is always reflection in the image side. The $k$-part of each transformation is understood through a six character code, which can be decoded
using the above table, showing what each characters $k$-part is.

We would like to point out that we have not said anything about the proof of why these side pairing transformations on the polyhedron $P$ actually
lead to hyperbolic 4-manifolds. This is all explained in their paper \cite{ratcliffe} from p.109-117, with details on the side-pairing coding from p.112-117.
We are taking it for granted that such manifolds exist and are well defined via side pairing transformations as described above.

Towards the end of their paper (see p. 117-124 \cite{ratcliffe}) Ratcliffe and Tschantz have tables of the 1171 manifolds giving various information about these 
manifolds, in particular the tables tell us what the side pairing code for each manifold is, and the link structure of each cusp. When we consider
explicit examples we will always write out explicitly what the code defining the manifold is and mention what the link structure of each cusp is.
The notation they use is to represent these link types are as follows. $A, B, \ldots , J$ represent the ten closed Euclidean 3-manifolds in the order given by 
Hantzsche and Wendt in  \cite{hantzsche}, these are also given in Wolf's text \cite{wolf} p.122, where
in his notation $\mathcal{G}_1$, \ldots , $\mathcal{G}_6$ correspond to $A, \ldots , F$ respectively, and $\mathcal{B}_1$, \ldots , $\mathcal{B}_4$ correspond
to $G, \ldots , J$. The orientable ones are given by $A, \ldots , F$ and the non-orientable ones are $G, \ldots , J$. Only $C, D, E$ do not occur as links
of cusps of any of their manifolds. Furthermore, $A$ is the 3-torus and $B$ is an orientable $S^1$-fibre bundle over the Klein bottle.

\section{A Handle Decomposition for the Ratcliffe-Tschantz manifolds}

In this section we want to explain how to obtain a handle decomposition for any one of the Ratcliffe-Tschantz manifolds. 

Recall from the previous section that each of the Ratcliffe-Tschantz manifolds is obtained via the 24-cell $P$ in $\mathbb{H}^4$ by appropriate side pairing 
transformations. Therefore it should come as no surprise that the approach to finding such a handle decomposition starts by trying to decompose the polyhedron
$P$ into various pieces that glue up appropriately when we apply the side pairing transformations. 
One must keep in mind that the polyhedron $P$ is not a manifold so we cannot really speak of a handle decomposition of it. The main idea is that
the dual cell decomposition we gave previously will lead to a handle decomposition for any of the Ratcliffe-Tschantz manifolds.

In section \ref{dual_24} we constructed a diagram of the dual polyhedron of $\partial{P}$. The red dots in the diagram will constitute pieces of the
1-handles for the manifolds. The lines joining these red dots will be parts of the 2-handles, and the triangles formed by these lines will
be parts of 3-handles. As for the 0-handle, just observe that the interior of $P$ can be identified with the interior of a 4-ball and as any
side pairing transformation is injective on the interior of $P$ we find that we can think of the interior as the 0-handle of any of the
Ratcliffe-Tschantz manifolds. Let us explain in slightly more detail how exactly we get this handle decomposition from the dual cell diagram
we drew.

Suppose we have obtained a manifold $M$ from side pairing transformations of the polyhedron $P$. The side pairing transformations form a group, which we denote
by $G$, we obtain $M$ by taking the quotient of $P$ by $G$. The polyhedron $P$ has a natural stratification into codimension
$k$ sides. The codimension 0 side is just the interior of $P$, so the zeroth level of the stratification of $P$ comprises of the interior of $P$. Codimension 1 sides
of $P$ comprise of the faces of $P$, and we know that there are 24 such faces, thus the first level of the stratification comprises of these 24 sides. There are
96 codimension 2 sides, hence the second level of the stratification comprises of these 96 codimension 2 sides. There are 96 codimension 3 sides, which 
tells us that the third level of the stratification consists of these 96 codimension 3 sides. The side pairing group $G$ preserves this stratification, and the orbit
of each codimension $k$ side in the $k^{th}$-level forms a codimension $k$ equivalence class. For example since $G$ does not identify any points in the
interior of $P$ we have that the codimension 0 equivalence class simply consists of the interior of $P$. Each side pairing transformation takes a side of $P$
and pairs it with another unique side of $P$, distinct from the original side. This implies that each codimension 1 equivalence class will consist of precisely
two elements, a pair of codimension 1 sides that are paired together by some side pairing transformation. For the codimension 2 sides there are exactly
four elements in each equivalence class (this is because the dihedral angles of a 24-cell are all $\pi/2$, hence four copies are needed around each codimension
$2$ side in order to give a total angle of $2\pi$), since there
are 96 codimension 2 sides it follows that there are a total of 24 codimension 2 equivalence classes. Finally, there are eight codimension 3 sides in each codimension
3 equivalence class, and since there are 96 codimension 3 sides there are a total of twelve codimension 3 equivalence classes.
From here it is easy to obtain a handle decomposition of $M$, each equivalence class of codimension $k$ sides
corresponds to one $k$ handle. As there is only one equivalence class for the codimension 0 side we have that $M$ has exactly
one 0-handle. As there are twelve equivalence classes of codimension 1 sides, it follows that $M$ has twelve 1-handles. Similarly one can conclude that $M$ has twenty four 
2-handles and twelve 3-handles. Finally, there are no 4-handles as all the vertices of $P$ are ideal.
We can now see why the dual cell diagram we drew in section \ref{dual_24} gives us the handle structure information. The vertices in the diagram corresponded to codimension 1 sides in $P$, hence pairs of vertices will correspond to 1-handles for $M$. Lines joining vertices corresponded to
codimension 2 sides, so collections of such lines (four to be exact) will form a 2-handle. Triangles filling in triples of lines that join vertices corresponded to
codimension 3 sides, so collections of such triangles (8 to be exact) will constitute a 3-handle.

It is time to start giving some explicit examples, so let us start with manifold no. 3 in their census. The code for this manifold is
\textbf{1477B8}, decoding this we find that the side pairing transformations are given by:

\[\xymatrixcolsep{5pc}\xymatrix{ S_{(+1,+1,0,0)}  \ar[r]^a_{k_{_{(-1,+1,+1,+1)}}} & S_{(-1,+1,0,0)} } \hspace{2cm} \xymatrix{S_{(+1,-1,0,0)}  \ar[r]^b_{k_{_{(-1,+1,+1,+1)}}} & S_{(-1,-1,0,0)} } \]

\[\xymatrixcolsep{5pc}\xymatrix{ S_{(+1,0,+1,0)}  \ar[r]^c_{k_{_{(+1,+1,-1,+1)}}} & S_{(+1,0,-1,0)} } \hspace{2cm} \xymatrix{S_{(-1,0,+1,0)}  \ar[r]^d_{k_{_{(+1,+1,-1,+1)}}} & S_{(-1,0,-1,0)} } \]

\[\xymatrixcolsep{5pc} \xymatrix{ S_{(0,+1,+1,0)}  \ar[r]^e_{k_{_{(-1,-1,-1,+1)}}} & S_{(0,-1,-1,0)} } \hspace{2cm} \xymatrix{S_{(0,+1,-1,0)}  \ar[r]^f_{k_{_{(-1,-1,-1,+1)}}} & S_{(0,-1,+1,0)} } \]

\[\xymatrixcolsep{5pc} \xymatrix{ S_{(+1,0,0,+1)}  \ar[r]^g_{k_{_{(-1,-1,-1,+1)}}} & S_{(-1,0,0,+1)} } \hspace{2cm} \xymatrix{S_{(+1,0,0,-1)}  \ar[r]^h_{k_{_{(-1,-1,-1,+1)}}} & S_{(-1,0,0,-1)} } \]

\[\xymatrixcolsep{5pc} \xymatrix{ S_{(0,+1,0,+1)}  \ar[r]^i_{k_{_{(-1,-1,+1,-1)}}} & S_{(0,-1,0,-1)} } \hspace{2cm} \xymatrix{S_{(0,+1,0,-1)}  \ar[r]^j_{k_{_{(-1,-1,+1,-1)}}} & S_{(0,-1,0,+1)} } \]

\[\xymatrixcolsep{5pc} \xymatrix{ S_{(0,0,+1,+1)}  \ar[r]^k_{k_{_{(+1,+1,+1,-1)}}} & S_{(0,0,+1,-1)} } \hspace{2cm} \xymatrix{S_{(0,0,-1,+1)}  \ar[r]^l_{k_{_{(+1,+1,+1,-1)}}} & S_{(0,0,-1,-1)} }. \]

From our above discussion we know that the 1-handles of this manifold correspond to pairs of identified sides of $P$. We are going to label the 1-handles as follows, 
each side pairing transformation pairs a domain side to an image side. We label the domain side by a capital letter, the letter corresponding to the letter of
the side pairing transformation. We label the image side by a primed capital letter, the letter being the same as what we used for the domain side. For example, 
from the above we see that the transformation $a$ pairs $S_{(+1,+1,0,0)}$ to $S_{(-1,+1,0,0)}$, hence we label $S_{(+1,+1,0,0)}$ by $A$ and
$S_{(-1,+1,0,0)}$ by $A'$. The following table summarises this information for manifold no. 3. \\

\begin{tabular}{|l|l|l||l|l|l|}
 $A$ &  $S_{(+1,+1,0,0)}$ & $(\frac{1}{\sqrt{2}}, \frac{1}{\sqrt{2}}, 0)$  &  $A'$ & $S_{(-1,+1,0,0)}$ & $(\frac{-1}{\sqrt{2}}, \frac{1}{\sqrt{2}}, 0)$  \\
 $B$ &  $S_{(+1,-1,0,0)}$ & $(\frac{1}{\sqrt{2}}, \frac{-1}{\sqrt{2}}, 0)$ &  $B'$ & $S_{(-1,-1,0,0)}$  & $(\frac{-1}{\sqrt{2}}, \frac{-1}{\sqrt{2}}, 0)$ \\
 $C$ &  $S_{(+1,0,+1,0)}$ &  $(\frac{1}{\sqrt{2}}, 0, \frac{1}{\sqrt{2}})$ &  $C'$ & $S_{(+1,0,-1,0)}$ & $(\frac{1}{\sqrt{2}}, 0, \frac{-1}{\sqrt{2}})$ \\
 $D$ &  $S_{(-1,0,+1,0)}$ &  $(\frac{-1}{\sqrt{2}}, 0, \frac{1}{\sqrt{2}})$ &  $D'$ & $S_{(-1,0,-1,0)}$ & $(\frac{-1}{\sqrt{2}}, 0, \frac{-1}{\sqrt{2}})$ \\ 
 $E$ &  $S_{(0,+1,+1,0)}$ &  $(0, \frac{1}{\sqrt{2}} ,\frac{1}{\sqrt{2}})$ &  $E'$ & $S_{(0,-1,-1,0)}$ &  $(0, \frac{-1}{\sqrt{2}} ,\frac{-1}{\sqrt{2}})$ \\ 
 $F$ &  $S_{(0,+1,-1,0)}$ &   $(0, \frac{1}{\sqrt{2}} ,\frac{-1}{\sqrt{2}})$ &  $F'$ & $S_{(0,-1,+1,0)}$ &  $(0, \frac{-1}{\sqrt{2}} ,\frac{1}{\sqrt{2}})$ \\ 
 $G$ &  $S_{(+1,0,0,+1)}$ &  $(1 + \sqrt{2}, 0, 0)$ &  $G'$ & $S_{(-1,0,0,-1)}$ & $(1 - \sqrt{2}, 0, 0)$ \\ 
 $H$ &  $S_{(+1,0,0,-1)}$ &  $(-1 + \sqrt{2}, 0, 0)$ &  $H'$ & $S_{(-1,0,0,+1)}$ & $(-1 - \sqrt{2}, 0, 0)$ \\
 $I$ &  $S_{(0,+1,0,+1)}$ &  $(0, 1 + \sqrt{2}, 0)$ &  $I'$ & $S_{(0,-1,0,+1)}$ &  $(0, -1 - \sqrt{2}, 0)$ \\
 $J$ &  $S_{(0,+1,0,-1)}$ &   $(0, -1 + \sqrt{2}, 0)$ &  $J'$ & $S_{(0,-1,0,-1)}$ &  $(0, 1 - \sqrt{2}, 0)$  \\
 $K$ &  $S_{(0,0,+1,+1)}$ &   $(0, 0, 1 + \sqrt{2})$ &  $K'$ & $S_{(0,0,+1,-1)}$ &  $(0, 0, -1 + \sqrt{2})$   \\
 $L$ &  $S_{(0,0,-1,+1)}$ &   $(0, 0, -1 - \sqrt{2})$  &  $L'$ & $S_{(0,0,-1,-1)}$ &  $(0, 0, 1 - \sqrt{2})$  
\end{tabular}  \\ \\

The following diagram shows the 1-handles.

\centerline {\graphicspath{ {example_manifold_3/} } \includegraphics[width=8cm, height=10cm]{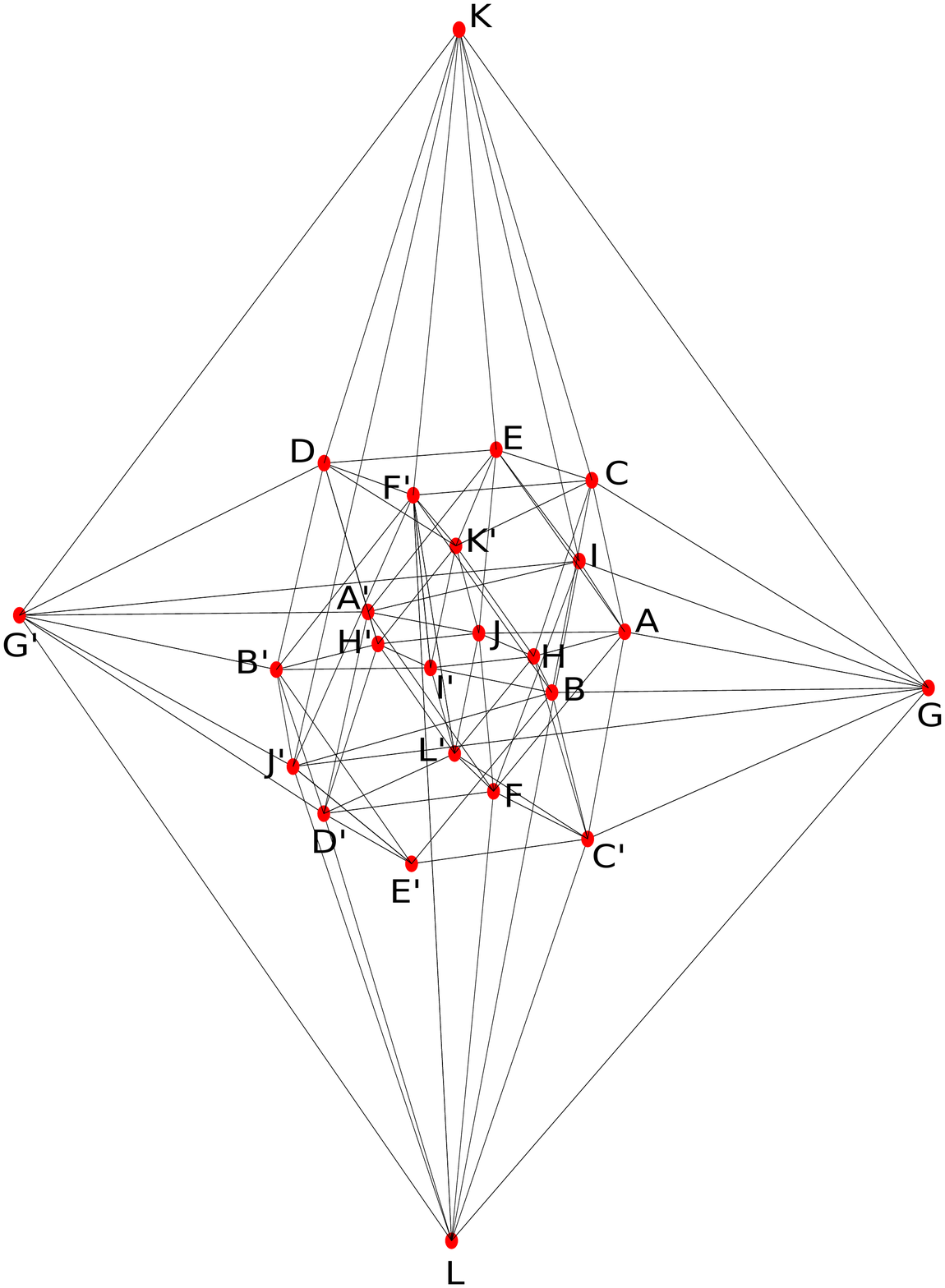}}

We should really think of the red vertices as very small 3-balls as a 1-handle is a copy of $D^1 \times D^3$, hence the attaching sphere is $S^0 \times D^3$, which
is a pair of 3-balls. The standard convention, for orientable manifolds, of visualising a 1-handle is by drawing its attaching region, which is $S^0 \times D^3$, in $\R^3$ with the
understanding that the boundaries of these two balls are identified via a reflection (an orientation reversing diffeomorphism). For the non-orientable case
the convention is that attaching regions of 1-handles are attached by an orientation preserving diffeomorphism. In our situation the 1-handles are in correspondence
with pairs of codimension one faces of $P$ (the 24-cell), and each such pair is identified via a side-pairing transformation. Therefore for our case the 
attaching regions of 1-handles will be identified via their corresponding side-pairing transformation, these could be orientation reversing
or preserving depending on the side pairing transformation. One has to be slightly careful in situations where the attaching regions are being identified via 
an orientation preserving diffeomorphism, the reason being that if we have a knotted 2-handle that passes over the 1-handle and we push it through the 1-handle, we find
that the knotted part of the 2-handle changes to its mirror. However, if all 2-handles present have components
that are unknotted then pushing any 2-handle through a 1-handle will not in anyway change the 2-handle. As we will see later we are in precisely
the situation where all the 2-handles are unknotted.
At this stage it is hard to get a good idea of how these 1-handles sit in 3-space, we will soon show a better way of visualising the 1-handles.

In order to see what a 2-handle looks like in the above picture we first need to understand how equivalence classes of codimension 2 sides arise. Recall that two codimension 1 sides intersect
if the centre vector of the spheres defining the sides have one entry the same, and the other entry in different positions. For example the side 
$S_{(+1,+1,0,0)}$ intersects $S_{(+1,0,+1,0)}$ because they both have a $+1$ in the first position and the second non-zero co-ordinates are in different positions.
On the other hand the side  $S_{(+1,+1,0,0)}$ does not intersect $S_{(0,0,-1,+1)}$ as their non-zero co-ordinates are in different positions. Also, 
$S_{(+1,+1,0,0)}$ does not intersect $S_{(-1,0,0,-1)}$ because the first entry of each have different signs.
Once we know which codimension 1 sides intersect, we find an equivalence class by simply applying side pairing transformations to the intersection. For example, we know that
$S_{(+1,+1,0,0)}$ intersects $S_{(+1,0,+1,0)}$, which using the above table we can write as $A \cap C$, we then apply the transformation $a$ and find that
$A \cap C$ goes to $A' \cap D$. We can then apply $d$ and we find that we end up at $A' \cap D'$. We continue in this way until we cycle back to $A \cap C$. Once we
end up back at $A \cap C$ we know we have found an equivalence class. The full equivalence class for $A \cap C$ is:

\[ \xymatrix{A \cap C \ar[r]^a & A'\cap D\ar[r]^d & A'\cap D' \ar[r]^{a^{-1}} & A\cap C' \ar[r]^{c^{-1}} & A\cap C } \]

Note that if $a$ pairs $A$ with $A'$ then $a^{-1}$ is the transformation that pairs $A'$ with $A$.
Continuing in this fashion we can work out all codimension 2 equivalence classes, the following table gives all 24 such classes.

\begin{tabular}{|l | l | }
\hline			
1. &  $\xymatrix{
A \cap C \ar[r]^a & A'\cap D\ar[r]^d & A'\cap D' \ar[r]^{a^{-1}} & A\cap C' \ar[r]^{c^{-1}} & A\cap C }$  \\ \hline
  
2. &  $\xymatrix{
A \cap G \ar[r]^a & A'\cap G' \ar[r]^{g^{-1}} & B\cap G \ar[r]^{b} & B'\cap G' \ar[r]^{g^{-1}} & A\cap G }$  \\ \hline  
  
3. &  $\xymatrix{
A \cap H \ar[r]^a & A'\cap H'\ar[r]^{h^{-1}} & B\cap H \ar[r]^{b} & B'\cap H' \ar[r]^{h^{-1}} & A\cap H }$  \\ \hline

4. &  $\xymatrix{
A \cap E \ar[r]^a & A'\cap E\ar[r]^e & B\cap E' \ar[r]^{b} & B'\cap E' \ar[r]^{e^{-1}} & A\cap E }$  \\ \hline

5. &  $\xymatrix{
A \cap F \ar[r]^a & A'\cap F\ar[r]^f & B\cap F' \ar[r]^{b} & B'\cap F' \ar[r]^{f^{-1}} & A\cap F }$  \\ \hline

6. &  $\xymatrix{
A \cap I \ar[r]^a & A'\cap I\ar[r]^i & B\cap I' \ar[r]^{b} & B'\cap I' \ar[r]^{i^{-1}} & A\cap I }$  \\ \hline  

7. &  $\xymatrix{
A \cap J \ar[r]^a & A'\cap J\ar[r]^j & B\cap J' \ar[r]^{b^{}} & B'\cap J' \ar[r]^{j^{-1}} & A\cap J }$  \\ \hline

8. &  $\xymatrix{
B \cap C \ar[r]^b & B'\cap D\ar[r]^d & B'\cap D' \ar[r]^{b^{-1}} & B\cap C' \ar[r]^{c^{-1}} & B\cap C }$  \\ \hline

9. &  $\xymatrix{
C \cap G \ar[r]^c & C'\cap G\ar[r]^g & D\cap G' \ar[r]^{d} & D'\cap G' \ar[r]^{g^{-1}} & C\cap G }$  \\ \hline

10. &  $\xymatrix{
C \cap H \ar[r]^c & C'\cap H \ar[r]^h & D\cap H' \ar[r]^{d} & D'\cap H' \ar[r]^{h^{-1}} & C\cap H }$  \\ \hline

11. &  $\xymatrix{
C \cap E \ar[r]^c & C'\cap F \ar[r]^f & D\cap F' \ar[r]^{d} & D'\cap E' \ar[r]^{e^{-1}} & C\cap E }$  \\ \hline

12. &  $\xymatrix{
C \cap F' \ar[r]^c & C'\cap E' \ar[r]^{e^{-1}} & D\cap E \ar[r]^{d} & D'\cap F \ar[r]^{f} & C\cap F' }$  \\   \hline

13. &  $\xymatrix{
C \cap K \ar[r]^c & C'\cap L \ar[r]^{l} & C'\cap L' \ar[r]^{c^{-1}} & C\cap K' \ar[r]^{k^{-1}} & C\cap K }$  \\   \hline

14. &  $\xymatrix{
D \cap K \ar[r]^d & D'\cap L \ar[r]^{l} & D'\cap L' \ar[r]^{d^{-1}} & D\cap K' \ar[r]^{k^{-1}} & D\cap K }$  \\   \hline

15. &  $\xymatrix{
G \cap I \ar[r]^g & G'\cap J' \ar[r]^{j^{-1}} & H\cap J \ar[r]^{h} & H'\cap I' \ar[r]^{i^{-1}} & G\cap I }$  \\   \hline

16. &  $\xymatrix{
G \cap J' \ar[r]^g & G'\cap I \ar[r]^{i} & H\cap I' \ar[r]^{h} & H'\cap J \ar[r]^{j} & G\cap J' }$  \\   \hline

17. &  $\xymatrix{
G \cap K \ar[r]^g & G'\cap L \ar[r]^{l} & H'\cap L' \ar[r]^{h^{-1}} & H\cap K' \ar[r]^{k^{-1}} & G\cap K }$  \\   \hline

18. &  $\xymatrix{
G \cap L \ar[r]^g & G'\cap K \ar[r]^{k} & H'\cap K' \ar[r]^{h^{-1}} & H\cap L' \ar[r]^{l^{-1}} & G\cap L }$  \\   \hline

19. &  $\xymatrix{
E \cap I \ar[r]^e & E'\cap J' \ar[r]^{j^{-1}} & F\cap J \ar[r]^{f} & F'\cap I' \ar[r]^{i^{-1}} & E\cap I }$  \\   \hline

20. &  $\xymatrix{
E \cap J \ar[r]^e & E'\cap I' \ar[r]^{i^{-1}} & F\cap I \ar[r]^{f} & F'\cap J' \ar[r]^{j^{-1}} & E\cap J }$  \\   \hline

21. &  $\xymatrix{
E \cap K \ar[r]^e & E'\cap L \ar[r]^{l} & E'\cap L' \ar[r]^{e^{-1}} & E\cap K' \ar[r]^{k^{-1}} & E\cap K }$  \\   \hline

22. &  $\xymatrix{
F \cap L \ar[r]^f & F'\cap K \ar[r]^{k} & F'\cap K' \ar[r]^{f^{-1}} & F\cap L' \ar[r]^{l^{-1}} & F\cap L }$  \\   \hline

23. &  $\xymatrix{
I \cap K \ar[r]^i & I'\cap K' \ar[r]^{k^{-1}} & J'\cap K \ar[r]^{j^{-1}} & J\cap K' \ar[r]^{k^{-1}} & I\cap K }$  \\   \hline

24. &  $\xymatrix{
I \cap L \ar[r]^i & I'\cap L' \ar[r]^{l^{-1}} & J'\cap L \ar[r]^{j^{-1}} & J\cap L' \ar[r]^{l^{-1}} & I\cap L }$  \\   \hline

\end{tabular} \\ \\

The following picture shows the $2$ handle given by the equivalence class

\[ \xymatrix{A \cap C \ar[r]^a & A'\cap D\ar[r]^d & A'\cap D' \ar[r]^{a^{-1}} & A\cap C' \ar[r]^{c^{-1}} & A\cap C } \]

The green edges and their blue vertices outline how the class looks in the whole handle decomposition.

\centerline {\graphicspath{ {example_manifold_3/} } \includegraphics[width=7cm, height=9cm]{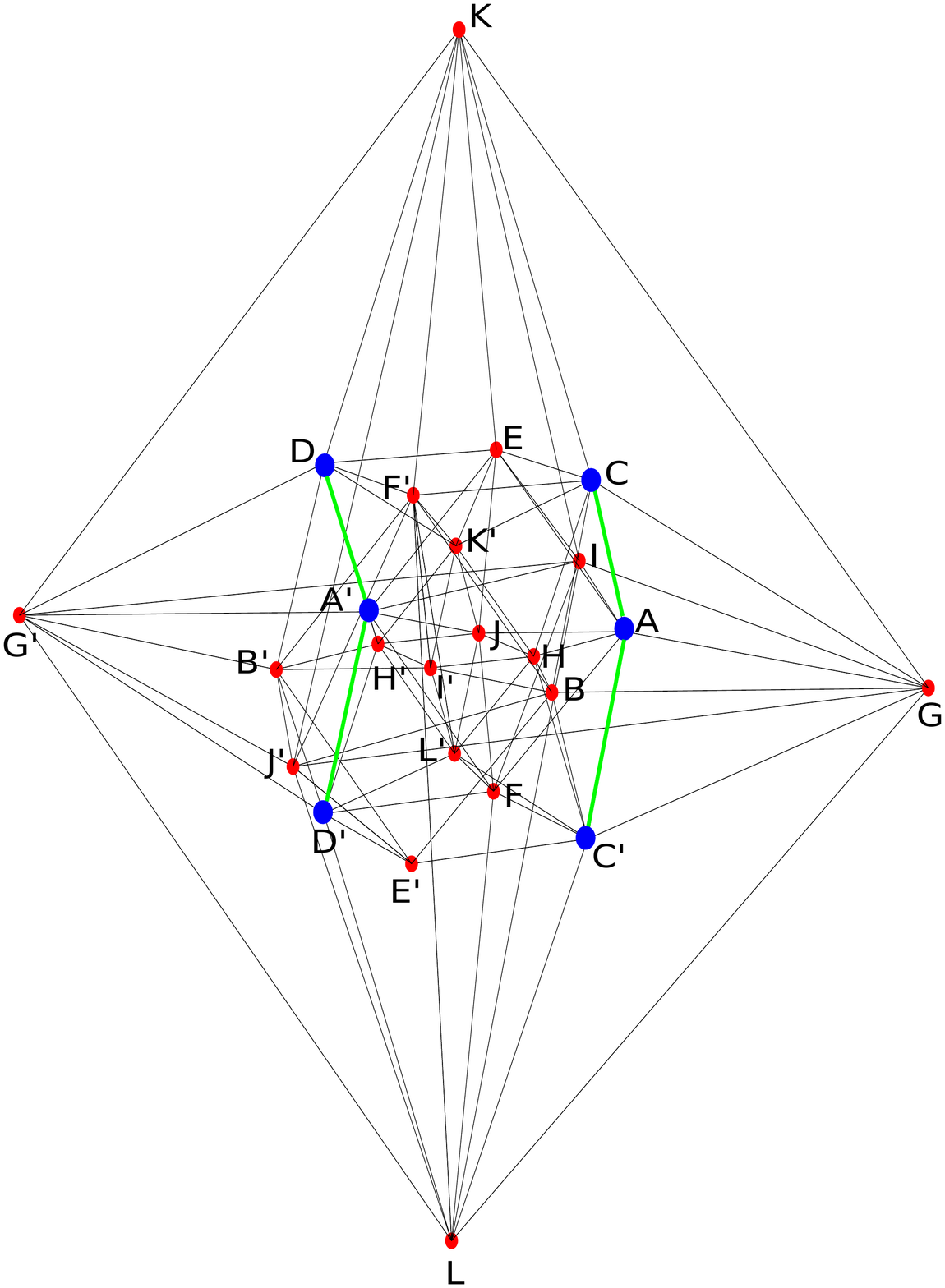}}

Let us give another picture of a $2$ handle, this time we take the 2-handle corresponding to the class
\[\xymatrix{
C \cap E \ar[r]^c & C'\cap F \ar[r]^f & D\cap F' \ar[r]^{d} & D'\cap E' \ar[r]^{e^{-1}} & C\cap E }\]

\centerline {\graphicspath{ {example_manifold_3/} } \includegraphics[width=7cm, height=9cm]{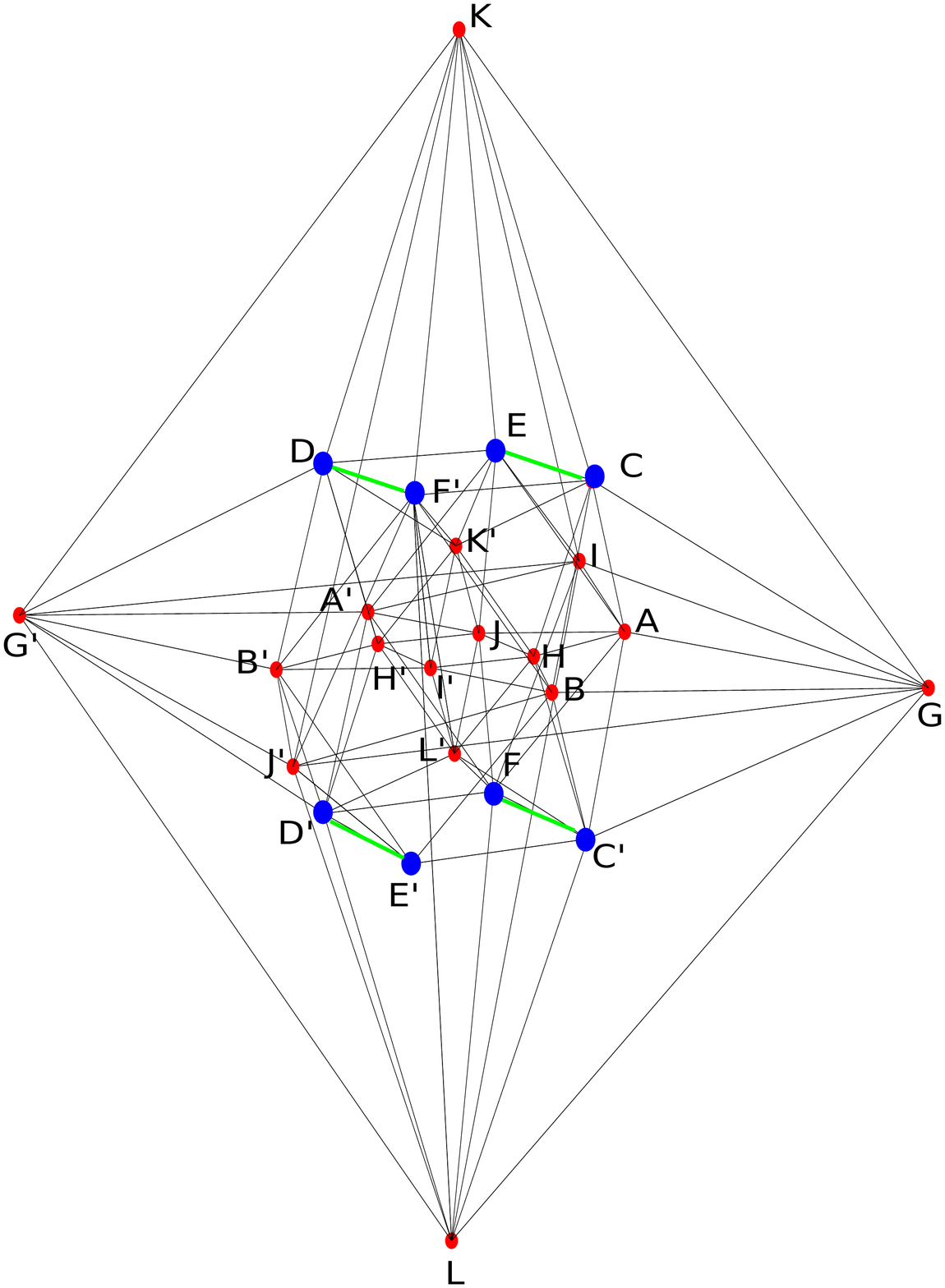}}

In total there are twenty four 2-handles and it is clear that drawing each 2-handle in a picture like the above would be very cumbersome. Also, drawing all of them
together in one picture proves to be a very difficult task, and even when one manages to do so it turns out to be very hard to see which bits of 2-handles belong
to which equivalence class.
In order to overcome these management issues of the 2-handles observe that if we split $\R^3$ into the $x-y$, $x-z$ and $y-z$ planes, then all of the 1-handles will lie
in these planes. For example the 1-handles $A-A'$ and $B-B'$ lie in the $x-y$ plane, $C-C'$ and $D-D'$ lie in the $x-z$ plane, $E-E'$ and $F-F'$ lie in the
$y-z$ plane. Some 1-handles will lie in the intersection of these planes, for example $I-I'$ lies in both the $x-y$ and $y-z$ planes. Exactly which plane (or planes) each
1-handle lies in can be found by simply looking at the table above showing the co-ordinates of the 1-handles. Due to the fact that each 1-handle lies in at least
one of three 2-planes, it turns out that many of the 2-handles will lie in one of these three 2-planes. The upshot of all of this is that if we split $\R^3$ into
these three 2-planes, we can then visualise the various 2-handles in each 2-plane separately. Unfortunately this hinges on the fact that every 2-handle resides on
one of these three 2-planes, and this is not true. It turns out that there are always six 2-handles which will not lie in any one of the $x-y$, $x-z$ or $y-z$ planes.
Therefore what we have to do is draw these six special 2-handles in a separate diagram, this means to show the structure of the 2-handles for each of the
Ratcliffe-Tschantz manifolds we will draw four diagrams. Three diagrams will consist of those 2-handles lying in the $x-y$, $x-z$, $y-z$ planes and one more diagram
showing the six special 2-handles that do not lie in any of these three 2-planes.

The following picture shows those 2-handles of manifold no. 3 that lie in the $x-y$ plane. We have drawn the 1-handles as 3-balls now as we have enough space to do so.

\centerline {\graphicspath{ {manifold_3/2-handle_cycles/} } \includegraphics[width=10cm, height=10cm]{2cycle_x-y}}

The 2-handles have been colour coded to make it easier to see which bits of 2-handle correspond to which 2-handle.

The following table shows which 2-handle corresponds to which colour:

\begin{tabular}{|l | l | }
\hline
colour & equivalence class  \\ \hline

green & $\xymatrix{
A \cap H \ar[r]^a & A'\cap H'\ar[r]^{h^{-1}} & B\cap H \ar[r]^{b} & B'\cap H' \ar[r]^{h^{-1}} & A\cap H }$ \\ \hline

red & $\xymatrix{
A \cap J \ar[r]^a & A'\cap J\ar[r]^j & B\cap J' \ar[r]^{b^{}} & B'\cap J' \ar[r]^{j^{-1}} & A\cap J }$  \\ \hline  

brown & $\xymatrix{
A \cap G \ar[r]^a & A'\cap G' \ar[r]^{g^{-1}} & B\cap G \ar[r]^{b} & B'\cap G' \ar[r]^{g^{-1}} & A\cap G }$   \\ \hline  

blue & $\xymatrix{
A \cap I \ar[r]^a & A'\cap I\ar[r]^i & B\cap I' \ar[r]^{b} & B'\cap I' \ar[r]^{i^{-1}} & A\cap I }$ \\ \hline

pink & $\xymatrix{
G \cap I \ar[r]^g & G'\cap J' \ar[r]^{j^{-1}} & H\cap J \ar[r]^{h} & H'\cap I' \ar[r]^{i^{-1}} & G\cap I }$  \\   \hline

black & $\xymatrix{
G \cap J' \ar[r]^g & G'\cap I \ar[r]^{i} & H\cap I' \ar[r]^{h} & H'\cap J \ar[r]^{j} & G\cap J' }$ \\   \hline

\end{tabular} \\ \\

The following picture shows those 2-handles that lie in the $x-z$ plane, with the table after it telling you which 2-handle corresponds to which colour.

\centerline {\graphicspath{ {manifold_3/2-handle_cycles/} } \includegraphics[width=10cm, height=10cm]{2cycle_x-z}}

\begin{tabular}{|l | l | }
\hline
colour & equivalence class  \\ \hline

green & $\xymatrix{
D \cap K \ar[r]^d & D'\cap L \ar[r]^{l} & D'\cap L' \ar[r]^{d^{-1}} & D\cap K' \ar[r]^{k^{-1}} & D\cap K }$  \\   \hline

red & $\xymatrix{
G \cap K \ar[r]^g & G'\cap L \ar[r]^{l} & H'\cap L' \ar[r]^{h^{-1}} & H\cap K' \ar[r]^{k^{-1}} & G\cap K }$  \\   \hline

brown & $\xymatrix{
C \cap G \ar[r]^c & C'\cap G\ar[r]^g & D\cap G' \ar[r]^{d} & D'\cap G' \ar[r]^{g^{-1}} & C\cap G }$  \\ \hline  

blue & $\xymatrix{
G \cap L \ar[r]^g & G'\cap K \ar[r]^{k} & H'\cap K' \ar[r]^{h^{-1}} & H\cap L' \ar[r]^{l^{-1}} & G\cap L }$  \\ \hline  

pink & $\xymatrix{
C \cap K \ar[r]^c & C'\cap L \ar[r]^{l} & C'\cap L' \ar[r]^{c^{-1}} & C\cap K' \ar[r]^{k^{-1}} & C\cap K }$  \\   \hline

black & $\xymatrix{
C \cap H \ar[r]^c & C'\cap H \ar[r]^h & D\cap H' \ar[r]^{d} & D'\cap H' \ar[r]^{h^{-1}} & C\cap H }$  \\   \hline

\end{tabular} \\ \\

The following picture shows those 2-handles that lie in the $y-z$ plane.

\centerline {\graphicspath{ {manifold_3/2-handle_cycles/} } \includegraphics[width=10cm, height=10cm]{2cycle_y-z}}

\begin{tabular}{|l | l | }
\hline
colour & equivalence class  \\ \hline

green & $\xymatrix{
I \cap L \ar[r]^i & I'\cap L' \ar[r]^{l^{-1}} & J'\cap L \ar[r]^{j^{-1}} & J\cap L' \ar[r]^{l^{-1}} & I\cap L }$  \\   \hline

red & $\xymatrix{
I \cap K \ar[r]^i & I'\cap K' \ar[r]^{k^{-1}} & J'\cap K \ar[r]^{j^{-1}} & J\cap K' \ar[r]^{k^{-1}} & I\cap K }$  \\   \hline

brown & $\xymatrix{
F \cap L \ar[r]^f & F'\cap K \ar[r]^{k} & F'\cap K' \ar[r]^{f^{-1}} & F\cap L' \ar[r]^{l^{-1}} & F\cap L }$ \\   \hline

blue & $\xymatrix{
E \cap K \ar[r]^e & E'\cap L \ar[r]^{l} & E'\cap L' \ar[r]^{e^{-1}} & E\cap K' \ar[r]^{k^{-1}} & E\cap K }$  \\   \hline

pink & $\xymatrix{
E \cap J \ar[r]^e & E'\cap I' \ar[r]^{i^{-1}} & F\cap I \ar[r]^{f} & F'\cap J' \ar[r]^{j^{-1}} & E\cap J }$  \\   \hline

black & $\xymatrix{
E \cap I \ar[r]^e & E'\cap J' \ar[r]^{j^{-1}} & F\cap J \ar[r]^{f} & F'\cap I' \ar[r]^{i^{-1}} & E\cap I }$  \\   \hline

\end{tabular} \\ \\

The final diagram shows a picture of those six 2-handles that do not lie in any one of the $x-y$, $x-z$, $y-z$ planes.

\centerline {\graphicspath{ {manifold_3/2-handle_cycles/} } \includegraphics[width=10cm, height=10cm]{2cycle_x-y-z}}

\begin{tabular}{|l | l | }
\hline
colour & equivalence class  \\ \hline

green & $\xymatrix{
A \cap E \ar[r]^a & A'\cap E\ar[r]^e & B\cap E' \ar[r]^{b} & B'\cap E' \ar[r]^{e^{-1}} & A\cap E }$  \\ \hline  

red & $\xymatrix{
B \cap C \ar[r]^b & B'\cap D\ar[r]^d & B'\cap D' \ar[r]^{b^{-1}} & B\cap C' \ar[r]^{c^{-1}} & B\cap C }$  \\ \hline  

brown & $\xymatrix{
A \cap C \ar[r]^a & A'\cap D\ar[r]^d & A'\cap D' \ar[r]^{a^{-1}} & A\cap C' \ar[r]^{c^{-1}} & A\cap C }$  \\ \hline

blue &  $\xymatrix{
A \cap F \ar[r]^a & A'\cap F\ar[r]^f & B\cap F' \ar[r]^{b} & B'\cap F' \ar[r]^{f^{-1}} & A\cap F }$  \\ \hline  

pink & $\xymatrix{
C \cap E \ar[r]^c & C'\cap F \ar[r]^f & D\cap F' \ar[r]^{d} & D'\cap E' \ar[r]^{e^{-1}} & C\cap E }$   \\ \hline  

black & $\xymatrix{
C \cap F' \ar[r]^c & C'\cap E' \ar[r]^{e^{-1}} & D\cap E \ar[r]^{d} & D'\cap F \ar[r]^{f} & C\cap F' }$ \\   \hline

\end{tabular} \\ \\

The above four diagrams showing the various 1 and 2-handles constitutes a complete Kirby diagram for manifold no. 3, and every Ratcliffe-Tschantz manifold
has a Kirby diagram that can be visualised by four similar diagrams. 

The 3-handles can be found in a similar way, the starting point is to work out all codimension 3 equivalence classes. A codimension 3 side is obtained by intersecting
three distinct codimension 1 sides. We already know that two codimension 1 sides intersect if the centre co-ordinates defining the associated spheres
have a common non-zero entry in the same position, and the other non-zero entry are in different positions. Three codimension 1 sides have non-empty intersection
if any pair of sides from the three sides have non-empty intersection. We then find an equivalence class as in the codimension 2 case, take such a non-empty intersection
and apply side pairing transformations until you cycle back to the original intersection of the three sides. Doing this for all codimension 3 sides
gives us all the equivalence classes, which in turn gives us all the 3-handles. 

The following table shows all the codimension 3 equivalence classes, there are twelve in total and each class contains eight sides.

\begin{table}[H]
  \resizebox{0.9\textwidth}{!}{\begin{minipage}{\textwidth}

\hskip-2.0cm\begin{tabular}{| l | l | }

\hline		
1. &  \tiny{$\xymatrix{
A \cap C \cap E \ar[r]^a & A'\cap D \cap E \ar[r]^e & B \cap C' \cap E' \ar[r]^{b} & B'\cap D' \cap E' \ar[r]^{d^{-1}} & B'\cap D \cap F'  \ar[r]^{b^{-1}} 
&  B\cap C \cap F' \ar[r]^{f^{-1}} & A' \cap D' \cap F \ar[r]^{a^{-1}} & A \cap C' \cap F}$}  \\ \hline
  
2. &  \tiny{$\xymatrix{
A \cap C \cap G \ar[r]^a & A'\cap D \cap G' \ar[r]^{g^{-1}} & B \cap C' \cap G \ar[r]^{b} & B'\cap D' \cap G' \ar[r]^{d^{-1}} & B'\cap D \cap G'  \ar[r]^{b^{-1}} 
&  B\cap C \cap G \ar[r]^{g} & A' \cap D' \cap G' \ar[r]^{a^{-1}} & A \cap C' \cap G}$}  \\ \hline
  
3. &  \tiny{$\xymatrix{
A \cap C \cap H \ar[r]^a & A'\cap D \cap H' \ar[r]^{h^{-1}} & B \cap C' \cap H \ar[r]^{b} & B'\cap D' \cap H' \ar[r]^{d^{-1}} & B'\cap D \cap H'  \ar[r]^{b^{-1}} 
&  B\cap C \cap H \ar[r]^{h} & A' \cap D' \cap H' \ar[r]^{a^{-1}} & A \cap C' \cap H}$}  \\ \hline

4.  &  \tiny{$\xymatrix{
A \cap E \cap I \ar[r]^a & A'\cap E \cap I \ar[r]^i & B \cap F' \cap I' \ar[r]^{b} & B'\cap F' \cap I' \ar[r]^{f^{-1}} & A\cap F \cap J  \ar[r]^{a} 
&  A'\cap F \cap J \ar[r]^{j} & B \cap E' \cap J' \ar[r]^b & B' \cap E' \cap J'}$}  \\ \hline

5. &  \tiny{$\xymatrix{
A \cap E \cap J \ar[r]^a & A'\cap E \cap J \ar[r]^j & B \cap F' \cap J' \ar[r]^{b} & B' \cap F' \cap J' \ar[r]^{f^{-1}} & A\cap F \cap I \ar[r]^{a} & A'\cap F \cap I  \ar[r]^{i} 
&  B\cap E' \cap I' \ar[r]^{b} & B' \cap E' \cap I'}$}  \\ \hline

6. &  \tiny{$\xymatrix{
A \cap G \cap I \ar[r]^a & A'\cap G' \cap I \ar[r]^i & B \cap H \cap I' \ar[r]^{b} & B'\cap H' \cap I' \ar[r]^{h^{-1}} & A\cap H \cap J  \ar[r]^{a} 
&  A'\cap H' \cap J \ar[r]^{j} & B \cap G \cap J' \ar[r]^b & B' \cap G' \cap J'}$}  \\ \hline

7. &  \tiny{$\xymatrix{
C \cap E \cap K \ar[r]^c & C'\cap F \cap L \ar[r]^l & C' \cap F \cap L' \ar[r]^{c^{-1}} & C\cap E \cap K' \ar[r]^{e} & D'\cap E' \cap L'  \ar[r]^{d^{-1}} 
&  D\cap F' \cap K' \ar[r]^{k^{-1}} & D \cap F' \cap K \ar[r]^d & D' \cap E' \cap L}$}  \\ \hline

8. &  \tiny{$\xymatrix{
C \cap F' \cap K \ar[r]^c & C'\cap E' \cap L \ar[r]^l & C' \cap E' \cap L' \ar[r]^{c^{-1}} & C\cap F' \cap K' \ar[r]^{f^{-1}} & D'\cap F \cap L'  \ar[r]^{d^{-1}} 
&  D\cap E \cap K' \ar[r]^{k^{-1}} & D \cap E \cap K \ar[r]^d & D' \cap F \cap L}$}  \\ \hline

9.  &  \tiny{$\xymatrix{
C \cap G \cap K \ar[r]^c & C'\cap G \cap L \ar[r]^l & C' \cap H \cap L' \ar[r]^{c^{-1}} & C\cap H \cap K' \ar[r]^{h} & D'\cap H' \cap L'  \ar[r]^{d^{-1}} 
&  D\cap H' \cap K' \ar[r]^{k^{-1}} & D \cap G' \cap K \ar[r]^d & D' \cap G' \cap L}$}  \\ \hline

10. &  \tiny{$\xymatrix{
E \cap I \cap K \ar[r]^e & E'\cap J' \cap L \ar[r]^{l} & E' \cap I' \cap L' \ar[r]^{e^{-1}} & E\cap J \cap K' \ar[r]^{j} & F'\cap J' \cap K  \ar[r]^{f^{-1}} 
&  F\cap I \cap L \ar[r]^{l} & F \cap J \cap L' \ar[r]^f & F' \cap I' \cap K'}$}  \\ \hline

11. &  \tiny{$\xymatrix{
G \cap I \cap K \ar[r]^g & G'\cap J' \cap L \ar[r]^{l} & H' \cap I' \cap L' \ar[r]^{h^{-1}} & H\cap J \cap K' \ar[r]^{j} & G'\cap J' \cap K  \ar[r]^{g^{-1}} 
&  G\cap I \cap L \ar[r]^{l} & H \cap J \cap L' \ar[r]^h & H' \cap I' \cap K'}$}  \\ \hline

12.  &  \tiny{$\xymatrix{
G \cap J' \cap K \ar[r]^g & G'\cap I \cap L \ar[r]^{l} & H' \cap J \cap L' \ar[r]^{h^{-1}} & H\cap I' \cap K' \ar[r]^{i^{-1}} & G'\cap I \cap K  \ar[r]^{g^{-1}} 
&  G\cap J' \cap L \ar[r]^{l} & H \cap I' \cap L' \ar[r]^h & H' \cap J \cap K'}$}  \\ \hline 

\end{tabular} 
      \end{minipage}}
\end{table}

As each class contains eight elements we find that each 3-handle will consist of eight triangles.

The following picture shows the 3-handle corresponding to the class:  
\[\hskip-1.5cm {\tiny{\xymatrix{
A \cap C \cap E \ar[r]^a & A'\cap D \cap E \ar[r]^e & B \cap C' \cap E' \ar[r]^{b} & B'\cap D' \cap E' \ar[r]^{d^{-1}} & B'\cap D \cap F'  \ar[r]^{b^{-1}} 
&  B\cap C \cap F' \ar[r]^{f^{-1}} & A' \cap D' \cap F \ar[r]^{a^{-1}} & A \cap C' \cap F}}} \]

\centerline {\graphicspath{ {manifold_3/3-handle_cycles/} } \includegraphics[width=10cm, height=10cm]{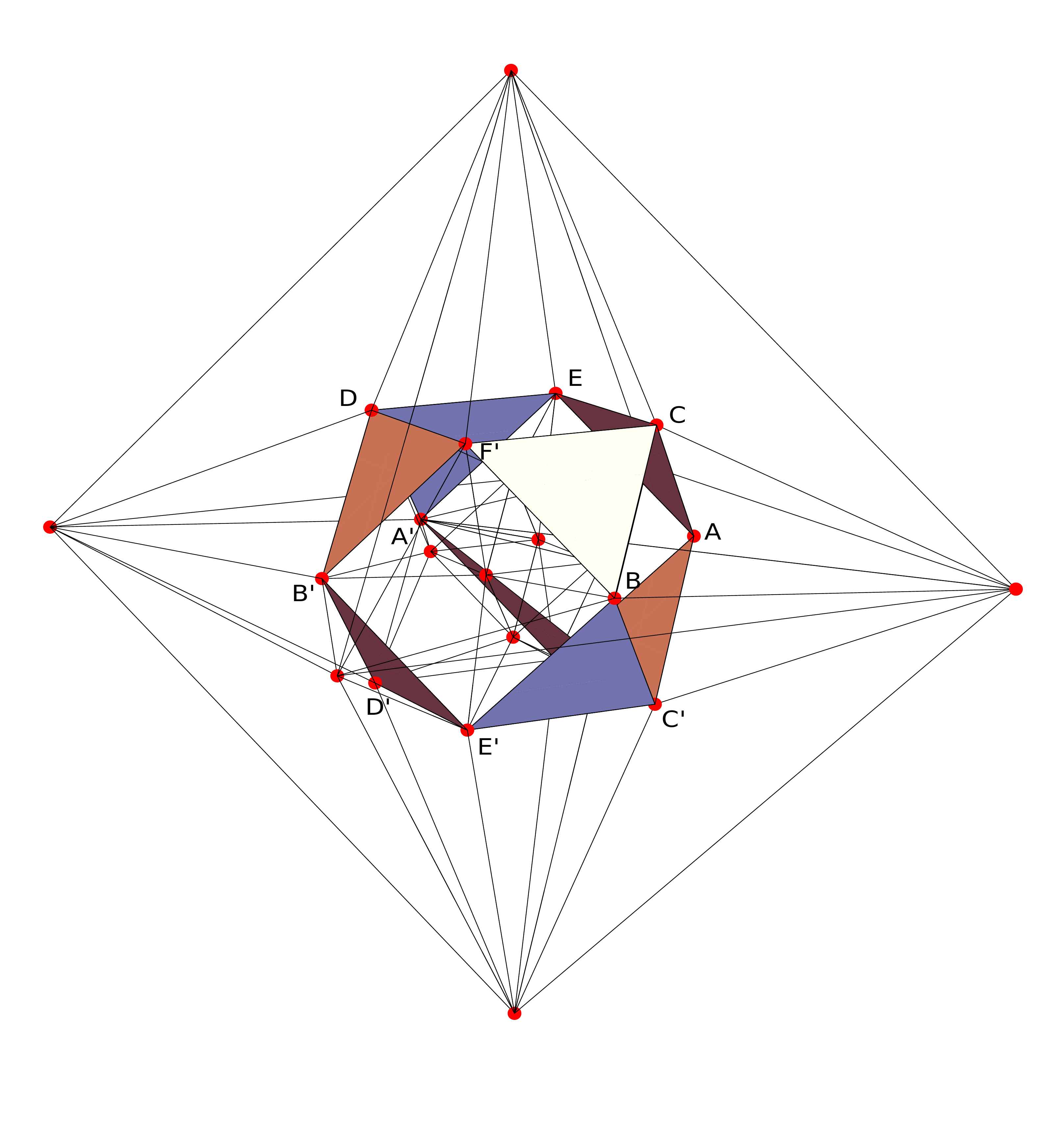}}

We have used a few colours to colour the triangles in the above picture to make it easier for the reader to see exactly where each triangle lies.
The totality of all eight triangles represents one 3-handle.

The following picture shows the 3-handle corresponding to the class:  
\[\hskip-1.5cm {\tiny{\xymatrix{
A \cap C \cap G \ar[r]^a & A'\cap D \cap G' \ar[r]^{g^{-1}} & B \cap C' \cap G \ar[r]^{b} & B'\cap D' \cap G' \ar[r]^{d^{-1}} & B'\cap D \cap G'  \ar[r]^{b^{-1}} 
&  B\cap C \cap G \ar[r]^{g} & A' \cap D' \cap G' \ar[r]^{a^{-1}} & A \cap C' \cap G}}} \]

\centerline {\graphicspath{ {manifold_3/3-handle_cycles/} } \includegraphics[width=10cm, height=10cm]{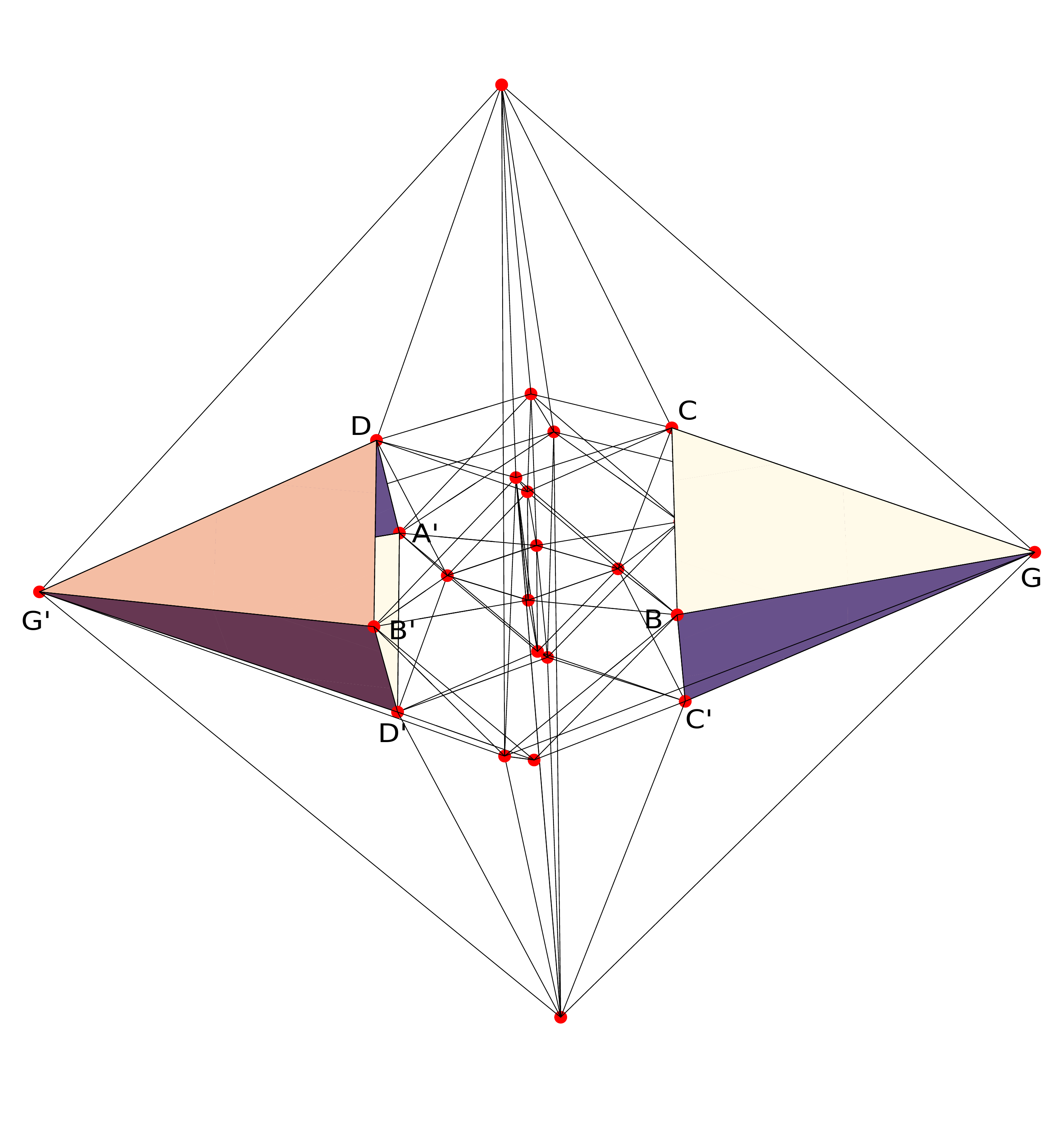}}

In the above picture we cannot see all components of the 3-handle, the bit of handle corresponding to $A \cap C \cap G$ lies behind $B \cap C \cap G$, and
the bit corresponding to $A \cap C' \cap G$ lies behind $B \cap C' \cap G$.

The following picture shows the 3-handle corresponding to the class:  
\[\hskip-2.5cm \tiny{\xymatrix{
A \cap C \cap H \ar[r]^a & A'\cap D \cap H' \ar[r]^{h^{-1}} & B \cap C' \cap H \ar[r]^{b} & B'\cap D' \cap H' \ar[r]^{d^{-1}} & B'\cap D \cap H'  \ar[r]^{b^{-1}} 
&  B\cap C \cap H \ar[r]^{h} & A' \cap D' \cap H' \ar[r]^{a^{-1}} & A \cap C' \cap H}} \]

\centerline {\graphicspath{ {manifold_3/3-handle_cycles/} } \includegraphics[width=9cm, height=9cm]{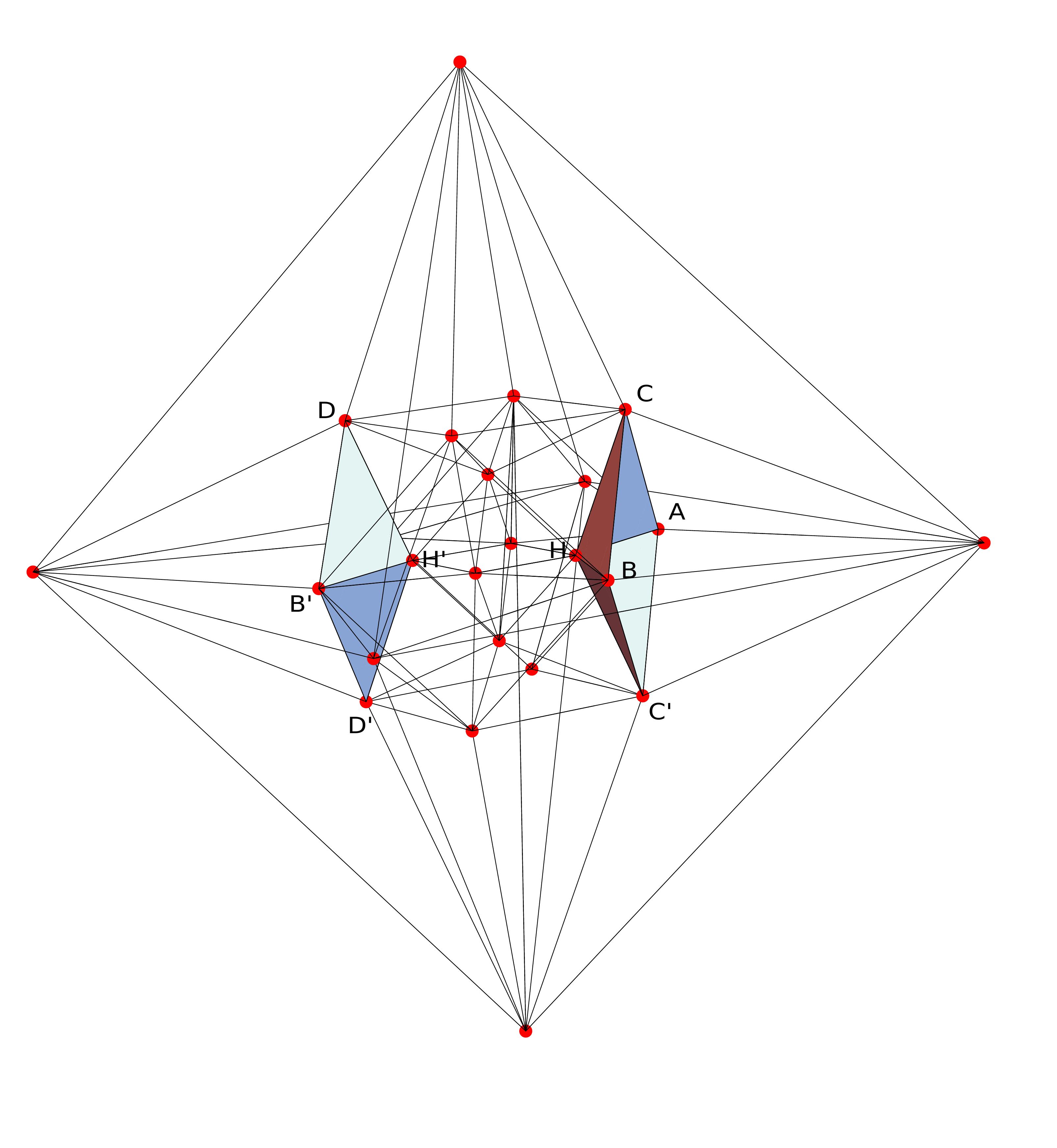}}

The triangles corresponding to $A' \cap D \cap H'$ and $A' \cap D' \cap H'$ lie behind $B' \cap D \cap H'$ and $B' \cap D' \cap H'$ respectively, and hence cannot
be seen in the above picture.
The above pictures of the first three 3-handles show that they do not lie in any one of the $x-y$, $x-z$ or $y-z$ planes. In general only three 3-handles will lie in any one of
the $x-y$, $x-z$ or $y-z$ planes. For manifold no. 3 the 3-handles corresponding to orbits of $A \cap G \cap I$, $C \cap G \cap K$ and $E \cap I \cap K$ 
(numbers 6, 9 and 10 in the above table) all lie in the $x-y$, $x-z$ and $y-z$ planes respectively

The following picture shows the 3-handle corresponding to the class determined by $A \cap G \cap I$.

\centerline {\graphicspath{ {manifold_3/3-handle_cycles/} } \includegraphics[width=8cm, height=8cm]{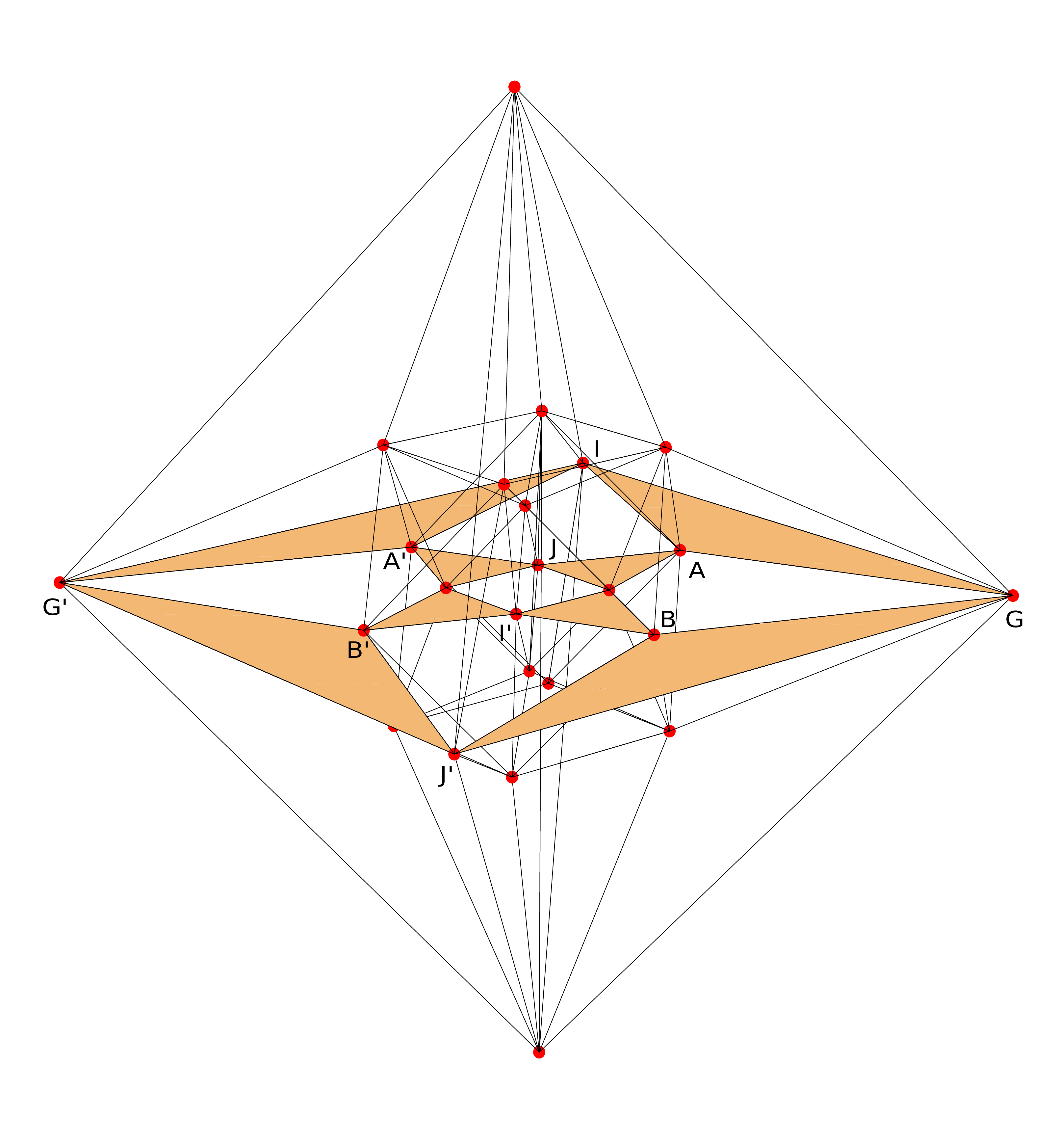}}

In summary we have shown the handle structure of the 4-manifold numbered 3 in the Ratcliffe-Tschantz census, furthermore
we have explicitly shown how the 2-handle structure can be clearly visualised via a Kirby diagram consisting of four diagrams.


Most of the Ratcliffe-Tschantz manifolds are non-orientable and one can ask whether
we can obtain Kirby diagrams in a similar fashion for all of them. The procedure is exactly analogous to the orientable case, for the sake of completeness
we give one non-orientable example.

Towards the beginning of their paper, Ratcliffe and Tschantz make special mention of one particular non-orientable manifold. They say 
``quite surprisingly, there is a congruence two 24-cell manifold with an even larger symmetry group of order 320. This manifold has the largest symmetry group among
all the congruence two 24-cell manifolds. If one equates beauty with symmetry, then this manifold is the most beautiful congruence two 24-cell manifold'' 
(see p. 102, second last paragraph in \cite{ratcliffe}). By ``congruence two 24-cell manifold'' they mean one of the Ratcliffe-Tschantz manifolds. 
The most beautiful manifold they speak of is numbered 1011 in the census, we believe such a comment gives enough incentive to investigate this manifold from the handle
decomposition viewpoint out of the possible 1149 choices. Therefore we start by describing, as in the case of manifold no. 3, its handle decomposition.

The side pairing code for manifold 1011 is \textbf{14FF28}, decoding this we obtain the explicit transformations:

\[\xymatrixcolsep{5pc} \xymatrix{ S_{(+1,+1,0,0)}  \ar[r]^a_{k_{(-1,+1,+1,+1)}} & S_{(-1,+1,0,0)} } \hspace{2cm} \xymatrix{S_{(+1,-1,0,0)}  \ar[r]^b_{k_{(-1,+1,+1,+1)}} & 
S_{(-1,-1,0,0)} }\]

\[\xymatrixcolsep{5pc} \xymatrix{ S_{(+1,0,+1,0)}  \ar[r]^c_{k_{(+1,+1,-1,+1)}} & S_{(+1,0,-1,0)} } \hspace{2cm} \xymatrix{S_{(-1,0,+1,0)}  \ar[r]^d_{k_{(+1,+1,-1,+1)}} & 
S_{(-1,0,-1,0)} }\]

\[\xymatrixcolsep{5pc} \xymatrix{ S_{(0,+1,+1,0)}  \ar[r]^e_{k_{(-1,-1,-1,-1)}} & S_{(0,-1,-1,0)} } \hspace{2cm} \xymatrix{S_{(0,+1,-1,0)}  \ar[r]^f_{k_{(-1,-1,-1,-1)}} & 
S_{(0,-1,+1,0)} }\]

\[\xymatrixcolsep{5pc} \xymatrix{ S_{(+1,0,0,+1)}  \ar[r]^g_{k_{(-1,-1,-1,-1)}} & S_{(-1,0,0,-1)} } \hspace{2cm} \xymatrix{S_{(+1,0,0,-1)}  \ar[r]^h_{k_{(-1,-1,-1,-1)}} & 
S_{(-1,0,0,+1)} }\]

\[\xymatrixcolsep{5pc} \xymatrix{ S_{(0,+1,0,+1)}  \ar[r]^i_{k_{(+1,-1,+1,+1)}} & S_{(0,-1,0,+1)} } \hspace{2cm} \xymatrix{S_{(0,+1,0,-1)}  \ar[r]^j_{k_{(+1,-1,+1,+1)}} & 
S_{(0,-1,0,-1)} }\]

\[\xymatrixcolsep{5pc} \xymatrix{ S_{(0,0,+1,+1)}  \ar[r]^k_{k_{(+1,+1,+1,-1)}} & S_{(0,0,+1,-1)} } \hspace{2cm} \xymatrix{S_{(0,0,-1,+1)}  \ar[r]^l_{k_{(+1,+1,+1,-1)}} & 
S_{(0,0,-1,-1)} }\]

The labelling of the 1-handles and their co-ordinates in $\R^3$ is shown in the following table: \\

\begin{tabular}{|l|l|l||l|l|l|}
 $A$ &  $S_{++00}$ & $(\frac{1}{\sqrt{2}}, \frac{1}{\sqrt{2}}, 0)$  &  $A'$ & $S_{-+00}$ & $(\frac{-1}{\sqrt{2}}, \frac{1}{\sqrt{2}}, 0)$  \\
 $B$ &  $S_{+-00}$ & $(\frac{1}{\sqrt{2}}, \frac{-1}{\sqrt{2}}, 0)$ &  $B'$ & $S_{--00}$  & $(\frac{-1}{\sqrt{2}}, \frac{-1}{\sqrt{2}}, 0)$ \\
 $C$ &  $S_{+0+0}$ &  $(\frac{1}{\sqrt{2}}, 0, \frac{1}{\sqrt{2}})$ &  $C'$ & $S_{+0-0}$ & $(\frac{1}{\sqrt{2}}, 0, \frac{-1}{\sqrt{2}})$ \\
 $D$ &  $S_{-0+0}$ &  $(\frac{-1}{\sqrt{2}}, 0, \frac{1}{\sqrt{2}})$ &  $D'$ & $S_{-0-0}$ & $(\frac{-1}{\sqrt{2}}, 0, \frac{-1}{\sqrt{2}})$ \\ 
 $E$ &  $S_{0++0}$ &  $(0, \frac{1}{\sqrt{2}} ,\frac{1}{\sqrt{2}})$ &  $E'$ & $S_{0--0}$ &  $(0, \frac{-1}{\sqrt{2}} ,\frac{-1}{\sqrt{2}})$ \\ 
 $F$ &  $S_{0+-0}$ &   $(0, \frac{1}{\sqrt{2}} ,\frac{-1}{\sqrt{2}})$ &  $F'$ & $S_{0-+0}$ &  $(0, \frac{-1}{\sqrt{2}} ,\frac{1}{\sqrt{2}})$ \\ 
 $G$ &  $S_{+00+}$ &  $(1 + \sqrt{2}, 0, 0)$ &  $G'$ & $S_{-00-}$ & $(1 - \sqrt{2}, 0, 0)$ \\ 
 $H$ &  $S_{+00-}$ &  $(-1 + \sqrt{2}, 0, 0)$ &  $H'$ & $S_{-00+}$ & $(-1 - \sqrt{2}, 0, 0)$ \\
 $I$ &  $S_{0+0+}$ &  $(0, 1 + \sqrt{2}, 0)$ &  $I'$ & $S_{0-0+}$ &  $(0, -1 - \sqrt{2}, 0)$ \\
 $J$ &  $S_{0+0-}$ &   $(0, -1 + \sqrt{2}, 0)$ &  $J'$ & $S_{0-0-}$ &  $(0, 1 - \sqrt{2}, 0)$  \\
 $K$ &  $S_{00++}$ &   $(0, 0, 1 + \sqrt{2})$ &  $K'$ & $S_{00+-}$ &  $(0, 0, -1 + \sqrt{2})$   \\
 $L$ &  $S_{00-+}$ &   $(0, 0, -1 - \sqrt{2})$  &  $L'$ & $S_{00--}$ &  $(0, 0, 1 - \sqrt{2})$  
\end{tabular} \\



The twenty four equivalence classes of $2$ handles is given in the following table with identifying transformations.

\begin{tabular}{|l | l | }
\hline			
1. &  $\xymatrix{
A \cap C \ar[r]^a & A'\cap D\ar[r]^d & A'\cap D' \ar[r]^{a^{-1}} & A\cap C' \ar[r]^{c^{-1}} & A\cap C }$  \\ \hline
  
2. &  $\xymatrix{
A \cap G \ar[r]^a & A'\cap H' \ar[r]^{h^{-1}} & B\cap H \ar[r]^{b} & B'\cap G' \ar[r]^{g^{-1}} & A\cap G }$  \\ \hline  
  
3. &  $\xymatrix{
A \cap H \ar[r]^a & A'\cap G'\ar[r]^{g^{-1}} & B\cap G \ar[r]^{b} & B'\cap H' \ar[r]^{h^{-1}} & A\cap H }$  \\ \hline

4. &  $\xymatrix{
A \cap E \ar[r]^a & A'\cap E\ar[r]^e & B\cap E' \ar[r]^{b} & B'\cap E' \ar[r]^{e^{-1}} & A\cap E }$  \\ \hline

5. &  $\xymatrix{
A \cap F \ar[r]^a & A'\cap F\ar[r]^f & B\cap F' \ar[r]^{b} & B'\cap F' \ar[r]^{f^{-1}} & A\cap F }$  \\ \hline

6. &  $\xymatrix{
A \cap I \ar[r]^a & A'\cap I\ar[r]^i & B'\cap i' \ar[r]^{b^{-1}} & B\cap I' \ar[r]^{i^{-1}} & A\cap I }$  \\ \hline  

7. &  $\xymatrix{
A \cap J \ar[r]^a & A'\cap J\ar[r]^j & B'\cap J' \ar[r]^{b^{-1}} & B\cap J' \ar[r]^{j^{-1}} & A\cap J }$  \\ \hline

8. &  $\xymatrix{
B \cap C \ar[r]^b & B'\cap D\ar[r]^d & B'\cap D' \ar[r]^{b^{-1}} & B\cap C' \ar[r]^{c^{-1}} & B\cap C }$  \\ \hline

9. &  $\xymatrix{
C \cap G \ar[r]^c & C'\cap G\ar[r]^g & D\cap G' \ar[r]^{d} & D'\cap G' \ar[r]^{g^{-1}} & C\cap G }$  \\ \hline

10. &  $\xymatrix{
C \cap H \ar[r]^c & C'\cap H \ar[r]^h & D\cap H' \ar[r]^{d} & D'\cap H' \ar[r]^{h^{-1}} & C\cap H }$  \\ \hline

11. &  $\xymatrix{
C \cap E \ar[r]^c & C'\cap F \ar[r]^f & D\cap F' \ar[r]^{d} & D'\cap E' \ar[r]^{e^{-1}} & C\cap E }$  \\ \hline

12. &  $\xymatrix{
C \cap F' \ar[r]^c & C'\cap E' \ar[r]^{e^{-1}} & D\cap E \ar[r]^{d} & D'\cap F \ar[r]^{f} & C\cap F' }$  \\   \hline

13. &  $\xymatrix{
C \cap K \ar[r]^c & C'\cap L \ar[r]^{l} & C'\cap L' \ar[r]^{c^{-1}} & C\cap K' \ar[r]^{k^{-1}} & C\cap K }$  \\   \hline

14. &  $\xymatrix{
D \cap K \ar[r]^d & D'\cap L \ar[r]^{l} & D'\cap L' \ar[r]^{d^{-1}} & D\cap K' \ar[r]^{k^{-1}} & D\cap K }$  \\   \hline

15. &  $\xymatrix{
G \cap I \ar[r]^g & G'\cap J' \ar[r]^{j^{-1}} & G'\cap J \ar[r]^{g^{-1}} & G\cap I' \ar[r]^{i^{-1}} & G\cap I }$  \\   \hline

16. &  $\xymatrix{
G \cap K \ar[r]^g & G'\cap L' \ar[r]^{l^{-1}} & H'\cap L \ar[r]^{h^{-1}} & H\cap K' \ar[r]^{k^{-1}} & G\cap K }$  \\   \hline

17. &  $\xymatrix{
G \cap L \ar[r]^g & G'\cap K' \ar[r]^{k^{-1}} & H'\cap K \ar[r]^{h^{-1}} & H\cap L' \ar[r]^{l^{-1}} & G\cap L }$  \\   \hline

18. &  $\xymatrix{
H \cap J \ar[r]^h & H'\cap I' \ar[r]^{i^{-1}} & H'\cap I \ar[r]^{h^{-1}} & H\cap J' \ar[r]^{j^{-1}} & H\cap J }$  \\   \hline

19. &  $\xymatrix{
E \cap I \ar[r]^e & E'\cap J' \ar[r]^{j^{-1}} & F\cap J \ar[r]^{f} & F'\cap I' \ar[r]^{i^{-1}} & E\cap I }$  \\   \hline

20. &  $\xymatrix{
E \cap J \ar[r]^e & E'\cap I' \ar[r]^{i^{-1}} & F\cap I \ar[r]^{f} & F'\cap J' \ar[r]^{j^{-1}} & E\cap J }$  \\   \hline

21. &  $\xymatrix{
E \cap K \ar[r]^e & E'\cap L' \ar[r]^{l^{-1}} & E'\cap L \ar[r]^{e^{-1}} & E\cap K' \ar[r]^{k^{-1}} & E\cap K }$  \\   \hline

22. &  $\xymatrix{
F \cap L \ar[r]^f & F'\cap K' \ar[r]^{k^{-1}} & F'\cap K \ar[r]^{f^{-1}} & F\cap L' \ar[r]^{l^{-1}} & F\cap L }$  \\   \hline

23. &  $\xymatrix{
I \cap K \ar[r]^i & I'\cap K \ar[r]^{k} & J'\cap K' \ar[r]^{j^{-1}} & J\cap K' \ar[r]^{k^{-1}} & I\cap K }$  \\   \hline

24. &  $\xymatrix{
I \cap L \ar[r]^i & I'\cap L \ar[r]^{l} & J'\cap L' \ar[r]^{j^{-1}} & J\cap L' \ar[r]^{l^{-1}} & I\cap L }$  \\   \hline

\end{tabular} \\ \\ \\

The following picture shows those 2-handles lying in the $x-y$ plane, with the table following explaining the colour coding.

\centerline {\graphicspath{ {2-handle_cycles/} }\includegraphics[width=10cm, height=10cm]{2cycle_x-y}}

\begin{tabular}{|l | l | }
\hline
colour & equivalence class  \\ \hline

green & $\xymatrix{
A \cap H \ar[r]^a & A'\cap G'\ar[r]^{g^{-1}} & B\cap G \ar[r]^{b} & B'\cap H' \ar[r]^{h^{-1}} & A\cap H }$ \\ \hline

red & $\xymatrix{
A \cap J \ar[r]^a & A'\cap J\ar[r]^j & B'\cap J' \ar[r]^{b^{-1}} & B\cap J' \ar[r]^{j^{-1}} & A\cap J }$  \\ \hline  

brown & $\xymatrix{
A \cap G \ar[r]^a & A'\cap H' \ar[r]^{h^{-1}} & B\cap H \ar[r]^{b} & B'\cap G' \ar[r]^{g^{-1}} & A\cap G }$  \\ \hline  

blue & $\xymatrix{
A \cap I \ar[r]^a & A'\cap I\ar[r]^i & B'\cap I' \ar[r]^{b^{-1}} & B\cap I' \ar[r]^{i^{-1}} & A\cap I }$ \\ \hline

pink & $\xymatrix{
G \cap I \ar[r]^g & G'\cap J' \ar[r]^{j^{-1}} & G'\cap J \ar[r]^{g^{-1}} & G\cap I' \ar[r]^{i^{-1}} & G\cap I }$  \\   \hline

black & $\xymatrix{
H \cap J \ar[r]^h & H'\cap I' \ar[r]^{i^{-1}} & H'\cap I \ar[r]^{h^{-1}} & H\cap J' \ar[r]^{j^{-1}} & H\cap J }$  \\   \hline

\end{tabular} \\ \\

The 2-handles that lie in the $x-z$ plane can be seen in the following picture.

\centerline {\graphicspath{ {2-handle_cycles/} }\includegraphics[width=10cm, height=10cm]{2cycle_x-z}}

\begin{tabular}{|l | l | }
\hline
colour & equivalence class  \\ \hline

green & $\xymatrix{
D \cap K \ar[r]^d & D'\cap L \ar[r]^{l} & D'\cap L' \ar[r]^{d^{-1}} & D\cap K' \ar[r]^{k^{-1}} & D\cap K }$  \\   \hline

red & $\xymatrix{
G \cap K \ar[r]^g & G'\cap L' \ar[r]^{l^{-1}} & H'\cap L \ar[r]^{h^{-1}} & H\cap K' \ar[r]^{k^{-1}} & G\cap K }$  \\   \hline

brown & $\xymatrix{
C \cap G \ar[r]^c & C'\cap G\ar[r]^g & D\cap G' \ar[r]^{d} & D'\cap G' \ar[r]^{g^{-1}} & C\cap G }$  \\ \hline  

blue & $\xymatrix{
C \cap H \ar[r]^c & C'\cap H \ar[r]^h & D\cap H' \ar[r]^{d} & D'\cap H' \ar[r]^{h^{-1}} & C\cap H }$  \\ \hline  

pink & $\xymatrix{
C \cap K \ar[r]^c & C'\cap L \ar[r]^{l} & C'\cap L' \ar[r]^{c^{-1}} & C\cap K' \ar[r]^{k^{-1}} & C\cap K }$  \\   \hline

black & $\xymatrix{
G \cap L \ar[r]^g & G'\cap K' \ar[r]^{k^{-1}} & H'\cap K \ar[r]^{h^{-1}} & H\cap L' \ar[r]^{l^{-1}} & G\cap L }$  \\   \hline

\end{tabular} \\ \\

The 2-handles that lie in the $y-z$ plane can be seen in the following picture.

\centerline {\graphicspath{ {2-handle_cycles/} }\includegraphics[width=10cm, height=10cm]{2cycle_y-z}}

\begin{tabular}{|l | l | }
\hline
colour & equivalence class  \\ \hline

green & $\xymatrix{
I \cap L \ar[r]^i & I'\cap L \ar[r]^{l} & J'\cap L' \ar[r]^{j^{-1}} & J\cap L' \ar[r]^{l^{-1}} & I\cap L }$  \\   \hline

red & $\xymatrix{
I \cap K \ar[r]^i & I'\cap K \ar[r]^{k} & J'\cap K' \ar[r]^{j^{-1}} & J\cap K' \ar[r]^{k^{-1}} & I\cap K }$  \\   \hline

brown & $\xymatrix{
F \cap L \ar[r]^f & F'\cap K' \ar[r]^{k^{-1}} & F'\cap K \ar[r]^{f^{-1}} & F\cap L' \ar[r]^{l^{-1}} & F\cap L }$  \\   \hline

blue & $\xymatrix{
E \cap K \ar[r]^e & E'\cap L' \ar[r]^{l^{-1}} & E'\cap L \ar[r]^{e^{-1}} & E\cap K' \ar[r]^{k^{-1}} & E\cap K }$  \\   \hline

pink & $\xymatrix{
E \cap J \ar[r]^e & E'\cap I' \ar[r]^{i^{-1}} & F\cap I \ar[r]^{f} & F'\cap J' \ar[r]^{j^{-1}} & E\cap J }$  \\   \hline

black & $\xymatrix{
E \cap I \ar[r]^e & E'\cap J' \ar[r]^{j^{-1}} & F\cap J \ar[r]^{f} & F'\cap I' \ar[r]^{i^{-1}} & E\cap I }$  \\   \hline

\end{tabular} \\ \\

Finally, the 2-handles that do not lie in any one of the above three planes can be seen in the following picture.

\centerline{\graphicspath{ {2-handle_cycles/} }\includegraphics[width=10cm, height=10cm]{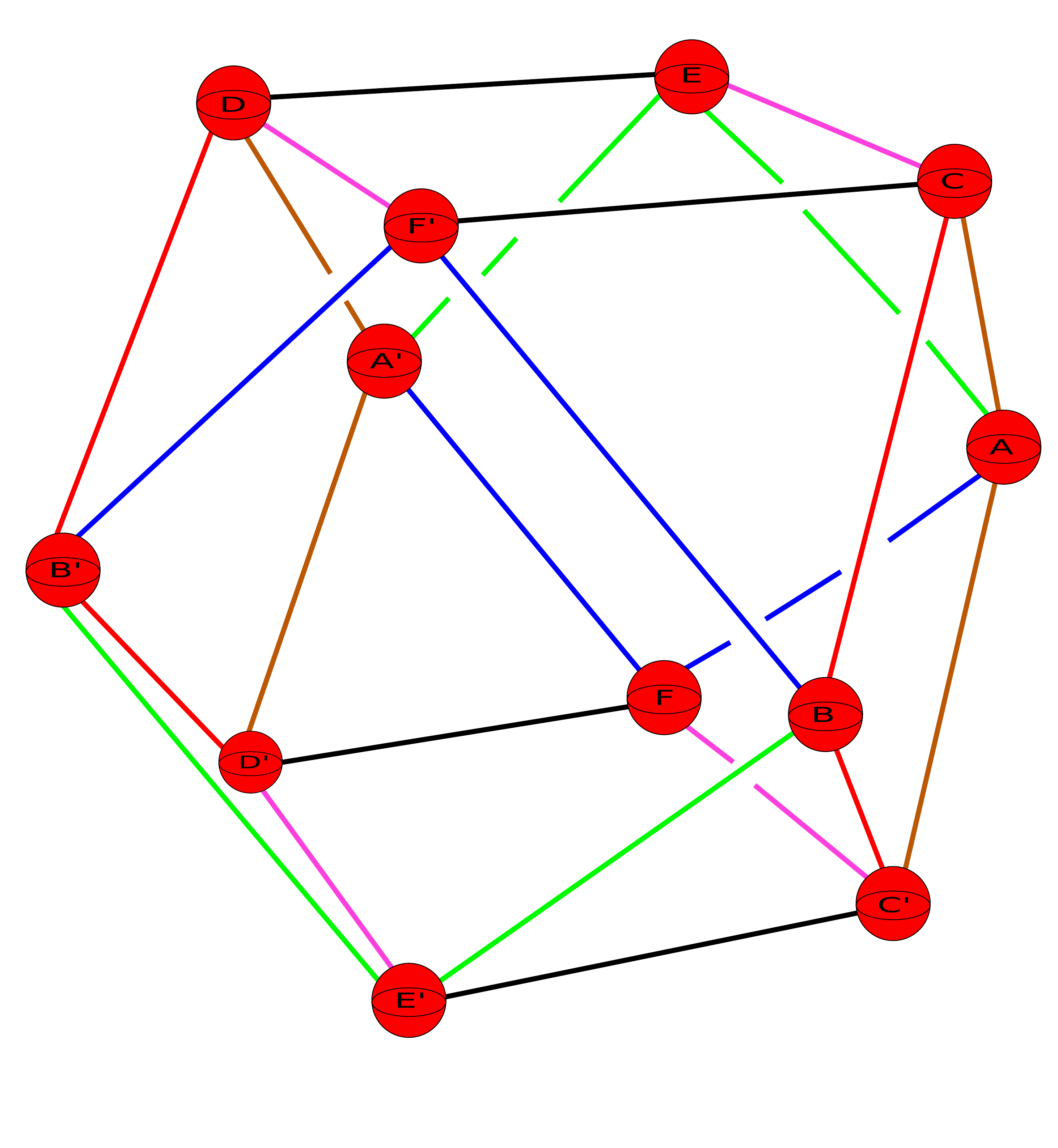}}

\begin{tabular}{|l | l | }
\hline
colour & equivalence class  \\ \hline

green & $\xymatrix{
A \cap E \ar[r]^a & A'\cap E\ar[r]^e & B\cap E' \ar[r]^{b} & B'\cap E' \ar[r]^{e^{-1}} & A\cap E }$  \\ \hline  

red & $\xymatrix{
B \cap C \ar[r]^b & B'\cap D\ar[r]^d & B'\cap D' \ar[r]^{b^{-1}} & B\cap C' \ar[r]^{c^{-1}} & B\cap C }$  \\ \hline  

brown & $\xymatrix{
A \cap C \ar[r]^a & A'\cap D\ar[r]^d & A'\cap D' \ar[r]^{a^{-1}} & A\cap C' \ar[r]^{c^{-1}} & A\cap C }$  \\ \hline

blue &  $\xymatrix{
A \cap F \ar[r]^a & A'\cap F\ar[r]^f & B\cap F' \ar[r]^{b} & B'\cap F' \ar[r]^{f^{-1}} & A\cap F }$  \\ \hline  

pink & $\xymatrix{
C \cap E \ar[r]^c & C'\cap F \ar[r]^f & D\cap F' \ar[r]^{d} & D'\cap E' \ar[r]^{e^{-1}} & C\cap E }$  \\ \hline  

black & $\xymatrix{
C \cap F' \ar[r]^c & C'\cap E' \ar[r]^{e^{-1}} & D\cap E \ar[r]^{d} & D'\cap F \ar[r]^{f} & C\cap F' }$  \\   \hline

\end{tabular} \\ \\

\section{The planar framing of the 2-handles}

A 2-handle in a 4-manifold is, by definition, simply a copy of $D^2 \times D^2$ and is attached along a copy of $S^1 \times D^2$. We can think
of the attaching data of a 2-handle to consist of two items:
\begin{itemize}

\item The attaching circle given by an embedding of $S^1 \times \{0\}$ in the 4-manifold, which is a copy of a knot in the one handle body structure.

\item A trivialisation of the normal bundle to the knot coming from the fact that we also have the $D^2$ factor to take in to account.

\end{itemize}
When drawing the attaching region of a 2-handle in a Kirby diagram, the general practise is to draw the attaching circle as some knot and then
give the information of the trivialisation of the normal bundle to the attaching circle in some external way. The trivialisation of the normal
bundle becomes important when one wants to carry out handle slides/cancellations. 
In the case that the attaching
circle of the 2-handle does not run over any 1-handles this is easy to do. Namely, the attaching circle is just a knot in some region of $\R^3$ and
the trivialisation of the associated normal bundle is given by an integer, called the \textit{framing number}. The reason why the trivialisation 
of the normal bundle corresponds to an integer is essentially because one can measure how twisted the normal bundle is by taking a parallel curve
to the attaching circle and counting (with sign) how many times it wraps around the attaching circle. One makes this intuition rigorous by
using the fact that a trivialisation of the normal bundle of the attaching circle of such a 2-handle corresponds to an element
of $\pi_1(O(2)) \cong \Z$, after one makes a choice of which element in  $\pi_1(O(2))$ corresponds to zero in $\Z$, for this
explanation see \cite{gompf} p.117. \\
In the general case, the attaching circle of a 2-handle may comprise of several components, consisting of arcs, with each component running over
a 1-handle. In this situation one has to be careful about how one encodes the trivialisation of the normal bundle. A natural instinct might
be to try and mimic the above case, where our 2-handle did not run across any 1-handles. However one must be very careful in trying to carry this process out,
in general such a generalisation of the \textit{framing number} is only well-defined modulo 2 (see \cite{gompf} p.120). 
There are quite a few ways to overcome these technicalities, we will not bother going in to the details of the various ways. Rather, we will
keep track of the trivialisation of the normal bundle to our 2-handles by defining what we call a \textit{planar framing}.

The following discussion will be valid for any one of the Ratcliffe-Tschantz manifolds. In two future papers we will be using manifolds 1011 and 35 together
with the constructions in this paper to construct examples of four dimensional hyperbolic link complements.
In view of this we will explain how to keep track of the normal bundle data 
associated to the attaching circles of the 2-handles of manifold 1011.

To start with we need to explicitly understand what the attaching maps for the attaching spheres of the 1-handles are. We are going to start by giving a description of what 
the attaching map for the 1-handle $A, A'$ looks like, this description will then provide
us with an explicit formula for computing the attaching map for all other 1-handles. The 1-handle $A, A'$ has attaching spheres consisting of one sphere 
centred at
$(\frac{1}{\sqrt{2}}, \frac{1}{\sqrt{2}}, 0)$, denoted by $A$ and another centred at $(\frac{-1}{\sqrt{2}}, \frac{1}{\sqrt{2}}, 0)$, denoted by $A'$.
Recall that the way we obtained the 1-handles in $\R^3$ was to identify each side of the 24-cell with a unique point on $S^3$ and then use stereographic
projection to map this point to a point in $\R^3$. Therefore, if we take a sphere of small radius about $(\frac{1}{\sqrt{2}}, \frac{1}{\sqrt{2}}, 0)$, applying
the inverse stereographic projection map we will obtain a small sphere centred about the point 
$(\frac{1}{\sqrt{2}}, \frac{1}{\sqrt{2}}, 0, 0)$ and contained in $S^3$. 
The line through $(0,0,0,0)$ and $(\frac{1}{\sqrt{2}}, \frac{1}{\sqrt{2}}, 0, 0)$ meets the sphere $S_{(+1,+1,0,0)}$ in a unique point, which we call $a_0$.

\centerline {\includegraphics[width=8cm, height=8cm]{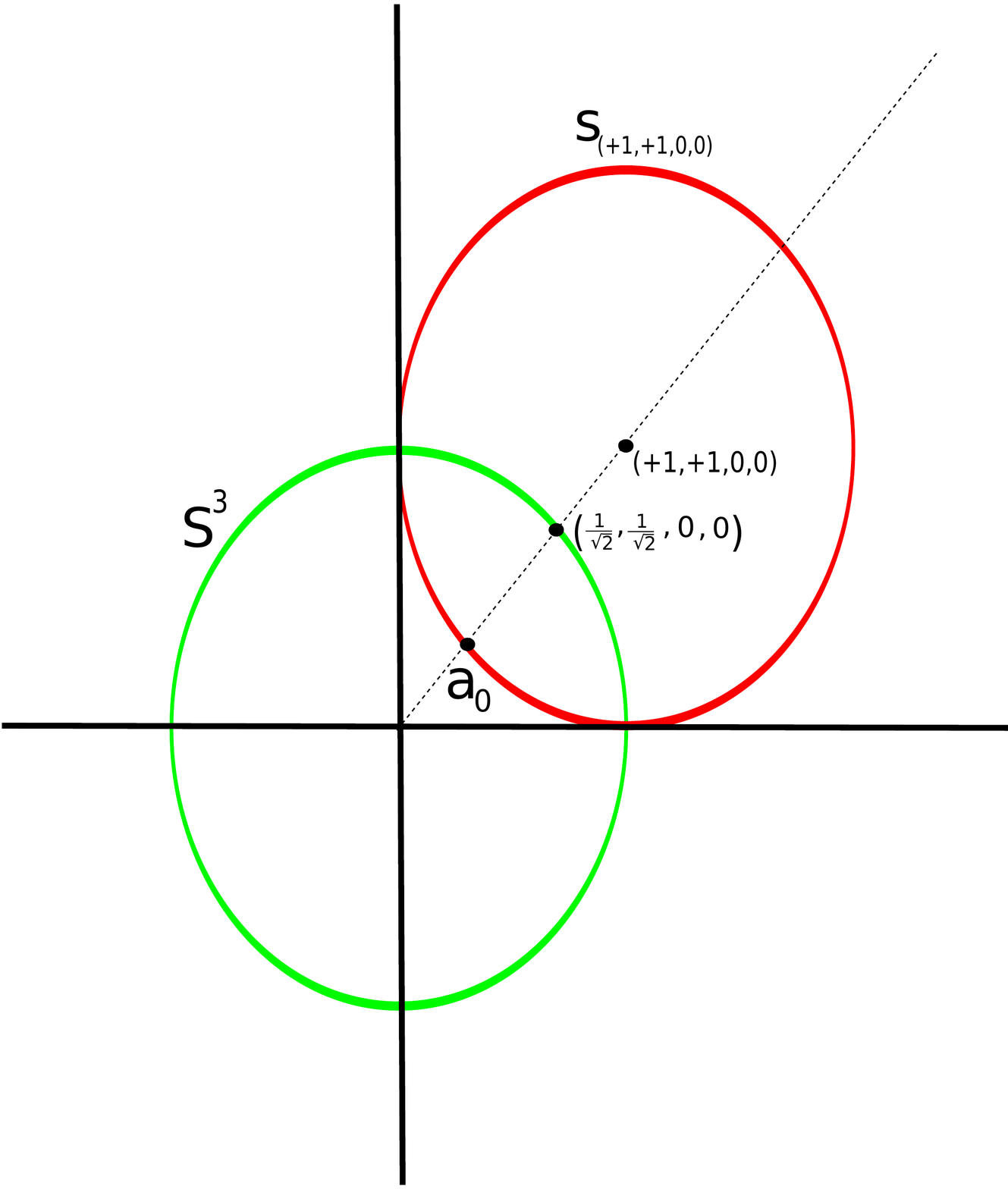}}

The sphere of small radius centred at $(\frac{1}{\sqrt{2}}, \frac{1}{\sqrt{2}}, 0, 0)$ can then be projected to a sphere of small radius about $a_0$.
When we apply the transformation $a = rk_{(-1,+1,+1,+1)}$, the small sphere about $a_0$ will get mapped to a small sphere about a point $a_1$ on
the side $S_{(-1,+1,0,0)}$. The point $a_1$ is given by the intersection of the line through $(0,0,0,0)$ and $(-\frac{1}{\sqrt{2}}, \frac{1}{\sqrt{2}}, 0, 0)$
with the sphere $S_{(-1,+1,0,0)}$. Observe that when applying $a$ the $r$-part is reflection in the image side, hence does not affect the small
sphere about $a_0$ and contained in $S_{(+1,+1,0,0)}$. We then radially project this small sphere about $a_1$ to obtain a small sphere about 
$(-\frac{1}{\sqrt{2}}, \frac{1}{\sqrt{2}}, 0, 0)$ and contained in $S^3$. Finally, applying the stereographic projection map we end up with
a small sphere about the point $(-\frac{1}{\sqrt{2}}, \frac{1}{\sqrt{2}}, 0)$.

In summary, we started with a small sphere centred about $(\frac{1}{\sqrt{2}}, \frac{1}{\sqrt{2}}, 0)$ in $\R^3$, corresponding to the side $A$, applied a diffeomorphism 
and obtained another sphere centred about $(-\frac{1}{\sqrt{2}}, \frac{1}{\sqrt{2}}, 0)$ in $\R^3$. This sphere corresponds to the side $A'$.
Note that the sphere centred at $(\frac{1}{\sqrt{2}}, \frac{1}{\sqrt{2}}, 0)$, corresponding to $A$, may have radius different from that
of the image sphere centred at $(-\frac{1}{\sqrt{2}}, \frac{1}{\sqrt{2}}, 0)$, corresponding to $A'$.

If we let $\phi : S^3 \rightarrow \R^3 \cup \{\infty\}$ denote stereographic projection, it is not hard to see that the above diffeomorphism
from the attaching sphere corresponding to $A$ to the attaching sphere corresponding to $A'$ is given by
\[\phi \circ k_{(-1,1,1,1)} \circ \phi^{-1}\]
This argument works for all 1-handles to show that the diffeomorphism that identifies the attaching spheres is given by
\[\phi \circ k \circ \phi^{-1}\]
where $k$ denotes the $k$-part of the transformation pairing the sides in question.

Using the above formula for the diffeomorphism we can explicitly work out how the attaching spheres of the 1-handles are identified.
We give the explicit details for the 1-handle $A, A'$. In this case the $k$-part is $k_{(-1,1,1,1)}$, we then need to work out $\phi \circ k_{(-1,1,1,1)} \circ \phi^{-1}$.
It is easy to compute this explicitly to obtain
\[ \phi \circ k_{(-1,1,1,1)} \circ \phi^{-1}(x,y,z) = (-x, y, z). \]
Therefore the attaching sphere corresponding to $A$, which is a sphere of small radius about $(\frac{1}{\sqrt{2}}, \frac{1}{\sqrt{2}}, 0)$ is identified
to the attaching sphere corresponding to $A'$, which is a sphere centred about $(-\frac{1}{\sqrt{2}}, \frac{1}{\sqrt{2}}, 0)$ via the reflection
$(x,y,z) \mapsto (-x,y,z)$.

We can carry out such computations for all 1-handles, the following table tells you exactly what the identifying diffeomorphism is in each case. \\

\begin{tabular}{|l|l|}
\hline
{} & {} \\
1-handle & Identifying diffeomorphism \\
{} & {} \\
\hline
{} & {} \\
$A,A'$ \& $B, B'$ & $(x,y,z) \mapsto (-x,y,z)$ \\
{} & {} \\
\hline
{} & {} \\
$C,C'$ \& $D, D'$ & $(x,y,z) \mapsto (x,y,-z)$ \\
{} & {} \\
\hline
{} & {} \\
$E,E'$ \& $F, F'$ & $(x,y,z) \mapsto \bigg(-\frac{x}{x^2+y^2+z^2},-\frac{y}{x^2+y^2+z^2},-\frac{z}{x^2+y^2+z^2}\bigg)$ \\
{} & {} \\
\hline
{} & {} \\
$G,G'$ \& $H, H'$ & $(x,y,z) \mapsto \bigg(-\frac{x}{x^2+y^2+z^2},-\frac{y}{x^2+y^2+z^2},-\frac{z}{x^2+y^2+z^2}\bigg)$ \\
{} & {} \\
\hline
{} & {} \\
$I,I'$ \& $J, J'$ & $(x,y,z) \mapsto (x,-y,z)$ \\
{} & {} \\
\hline
{} & {} \\
$K,K'$ \& $L, L'$ & $(x,y,z) \mapsto \bigg(\frac{x}{x^2+y^2+z^2},\frac{y}{x^2+y^2+z^2},\frac{z}{x^2+y^2+z^2}\bigg)$ \\
{} & {} \\
\hline

\end{tabular} \\

From the table you can see that all attaching spheres of 1-handles are identified by either a reflection or by a composition of a collection
of reflections together with an inversion in $S^2$.

Now that we know exactly how the attaching spheres of each 1-handle are identified we can deal with the problem of understanding the normal
bundle data associated to each 2-handle. We can split the 2-handles in to two groups, those that lie in the $x-y$, $x-z$ or $y-z$ planes, and 
the other six that do not all lie in any one of these planes. When carrying out a handle slide of a 2-handle over another 2-handle one has
to pick a parallel curve to the 2-handle, the handle slide is then done by carrying out a type of band sum (see \cite{gompf} p.139 for details).
It is in this regard that the trivialisation of the normal bundle of a 2-handle component becomes important as it tells us how the parallel curve
looks in our diagram. Therefore we need to consider two situations, how a trivialisation for the normal bundle of those 2-handles lying in the
$x-y$, $x-z$ or $y-z$ planes looks, and then how a trivialisation for the normal bundle for the six 2-handles that do not lie in any one of these
planes looks.

We will give the details for the six special 2-handles that do not lie in any one of the $x-y$, $x-z$ and $y-z$ planes. The same argument can be used
to show that any 2-handle lying in the $x-y$, $x-z$ or $y-z$ plane has a parallel curve that is confined to such a plane. The rough idea is as follows:
Start with the $x-y$ plane, each 2-handle has four components consisting of straight line segments. Fix any 2-handle in the $x-y$ plane, we can then take
a parallel straight line segment to one component of the 2-handle and chase it around. What we find is that because the attaching circle of this 2-handle is confined to
the $x-y$ plane any parallel curve will also be confined to the $x-y$ plane. What this tells us is that the trivialisation of the normal
bundle of a component of such a 2-handle is one that does not ``twist'' out of the $x-y$ plane i.e. it is confined to the $x-y$ plane.
A similar argument shows that all the 2-handles that reside in one of the $x-y$, $x-z$ or $y-z$ planes have components whose normal bundles are
trivialised in a similar manner, in particular this implies that any such parallel curve to such a 2-handle must remain confined to such a 
plane, and hence cannot twist around any component of the 2-handle in any way. It thus follows that if we were to carry out a handle slide
within one of these planes then the resulting 2-handle will also be confined to this plane.

The details for the six special 2-handles are as follows. Start with the diagram showing the six 2-handles that do not all lie in the $x-y$, $x-z$ or $y-z$ planes. What we 
need to do is understand how
a parallel curve behaves to each such 2-handle. It is here that it becomes vital that we explicitly know how the attaching spheres associated to each
1-handle are being identified. This is because we are going to take a parallel segment and then follow it along the 2-handle keeping track of how it emerges
out of attaching spheres of 1-handles.
Let us give the details in the case of the pink 2-handle. We start by taking a parallel curve to the pink 2-handle at the point shown in the following diagram.
We take a parallel curve just below the pink 2-handle component running from $C$ to $E$, shown as a black dashed line.

\centerline {\graphicspath{ {framing/} } \includegraphics[width=8cm, height=7cm]{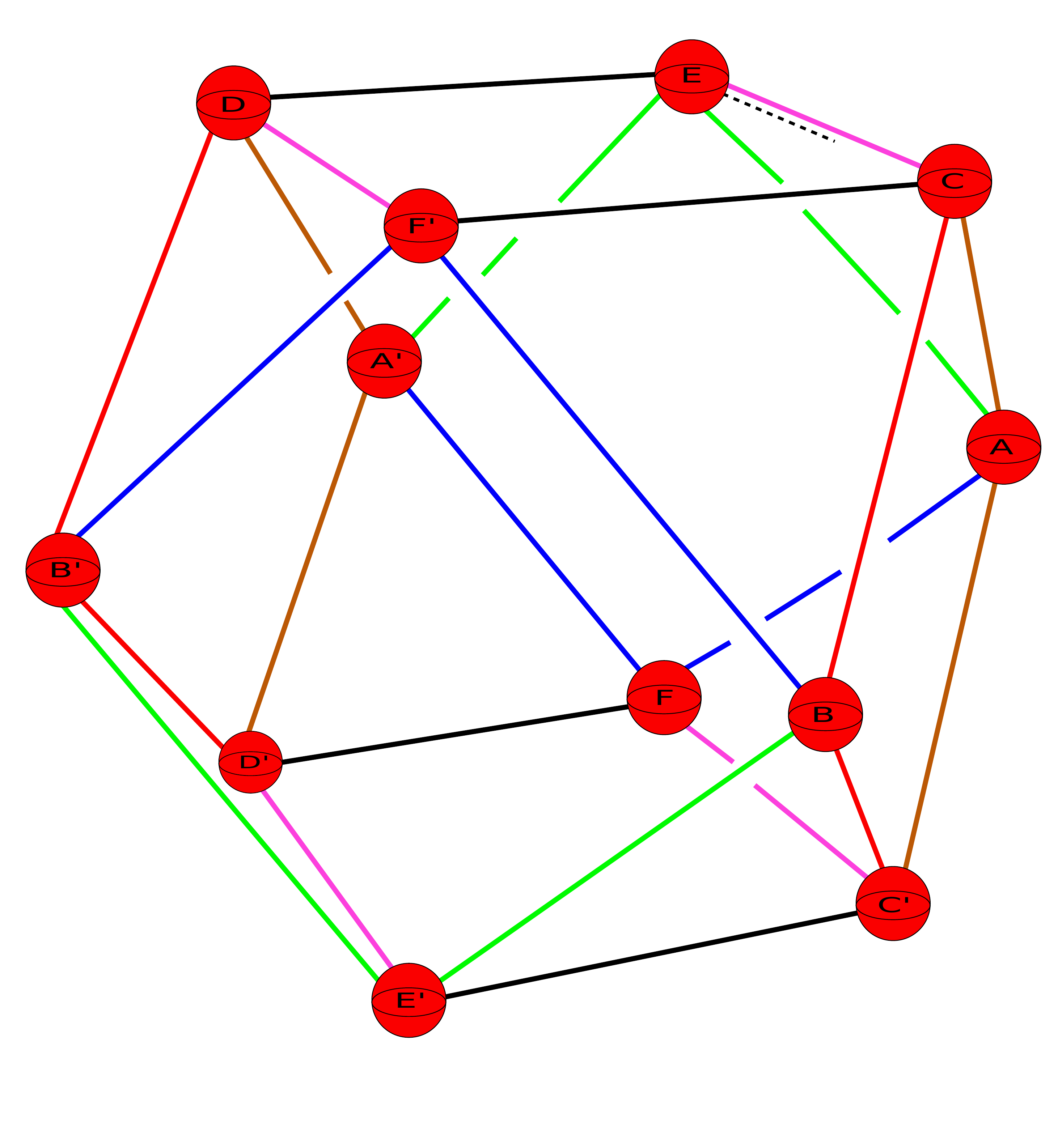}}

When this parallel curve hits $E$ it comes back out of $E'$, we need to work out exactly where it comes out of $E'$. Recall from the above we know
that the diffeomorphism that identifies the attaching sphere of $E$ to that of $E'$ is given by 
\[ (x,y,z) \mapsto \bigg(-\frac{x}{x^2+y^2+z^2},-\frac{y}{x^2+y^2+z^2},-\frac{z}{x^2+y^2+z^2}\bigg) \]
which is a composition of an inversion in the unit sphere centred at the origin together with the antipodal map. 
First of all, the pink 2-handle component that goes into $E$ comes out of $E'$ in the following way.

\centerline {\graphicspath{ {framing/} } \includegraphics[width=4cm, height=5cm]{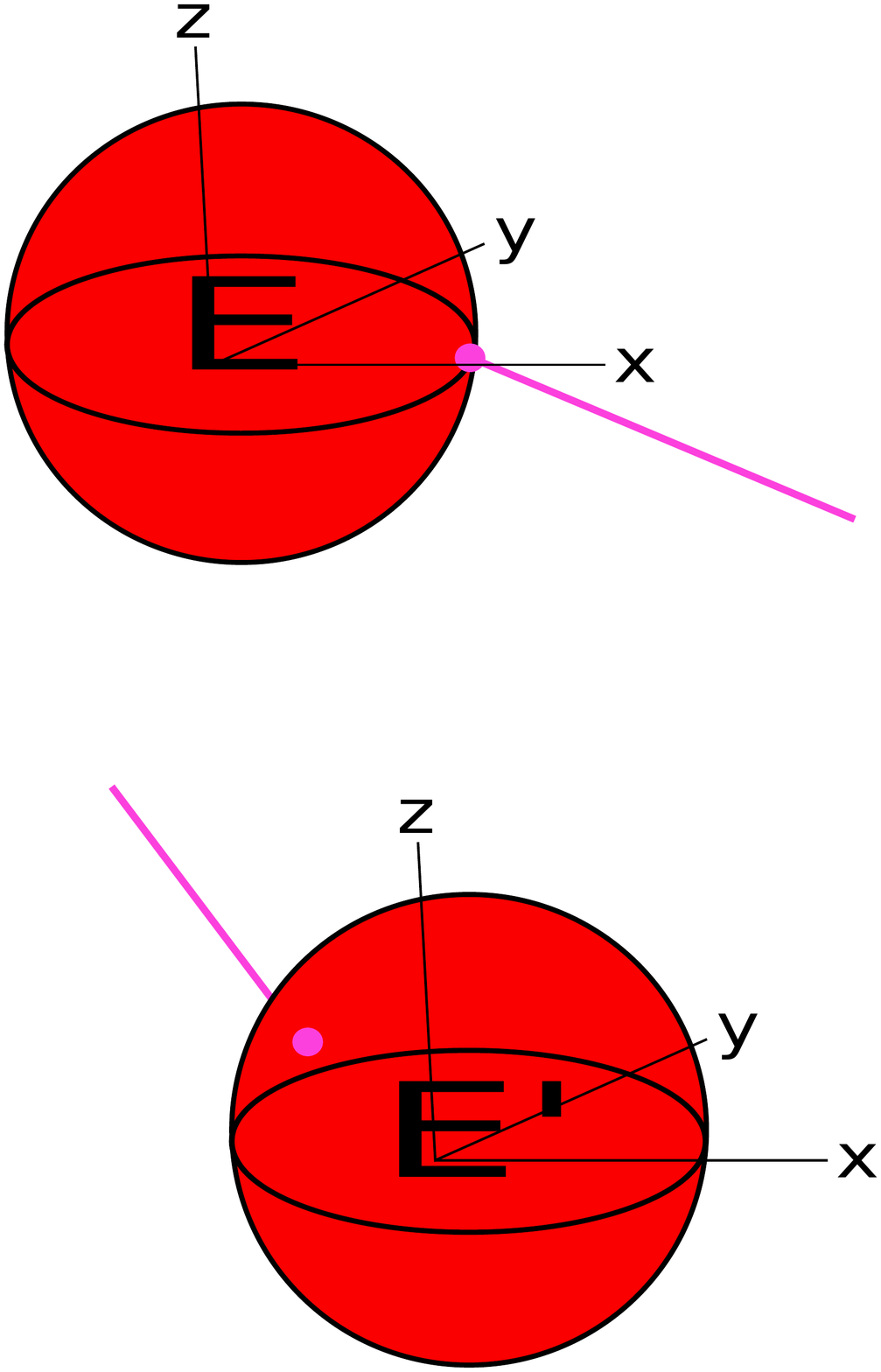}}

We have put axes at each centre point of each sphere so that the reader knows which directions are positive. To see how we get
this diagram, just observe that the point at which the pink 2-handle component hits $E$ can be written as
$(r, \frac{1}{\sqrt{2}}, \frac{1}{\sqrt{2}})$, where $r$ is the radius of the attaching sphere (some small number, which we do not need to explicitly know).
This point lies outside of the unit sphere $S^2 \subseteq \R^3$. If we apply the inversion map 
\[ (x,y,z) \mapsto (\frac{x}{x^2 + y^2 + z^2}, \frac{y}{x^2 + y^2 + z^2}, \frac{z}{x^2 + y^2 + z^2}) \]
The point $(r, \frac{1}{\sqrt{2}}, \frac{1}{\sqrt{2}})$ will get mapped to $(\frac{r}{r^2 + 1},\frac{1}{\sqrt{2}(r^2 + 1)}, \frac{1}{\sqrt{2}(r^2 + 1)})$, which
lies inside the unit sphere $S^2$. 
In fact, the attaching sphere $E$ centred at $(0, \frac{1}{\sqrt{2}}, \frac{1}{\sqrt{2}})$ will map to a sphere under the inversion map. One can carry out
some simple computations to work out that the image sphere will have centre $(0, \frac{1}{\sqrt{2}(1-r^2)}, \frac{1}{\sqrt{2}(1-r^2)})$ and have radius
$R := \frac{r}{1-r^2}$.
The following diagram shows the image of $(r, \frac{1}{\sqrt{2}}, \frac{1}{\sqrt{2}})$ under the inversion map.

\centerline {\graphicspath{ {framing/} } \includegraphics[width=6cm, height=6cm]{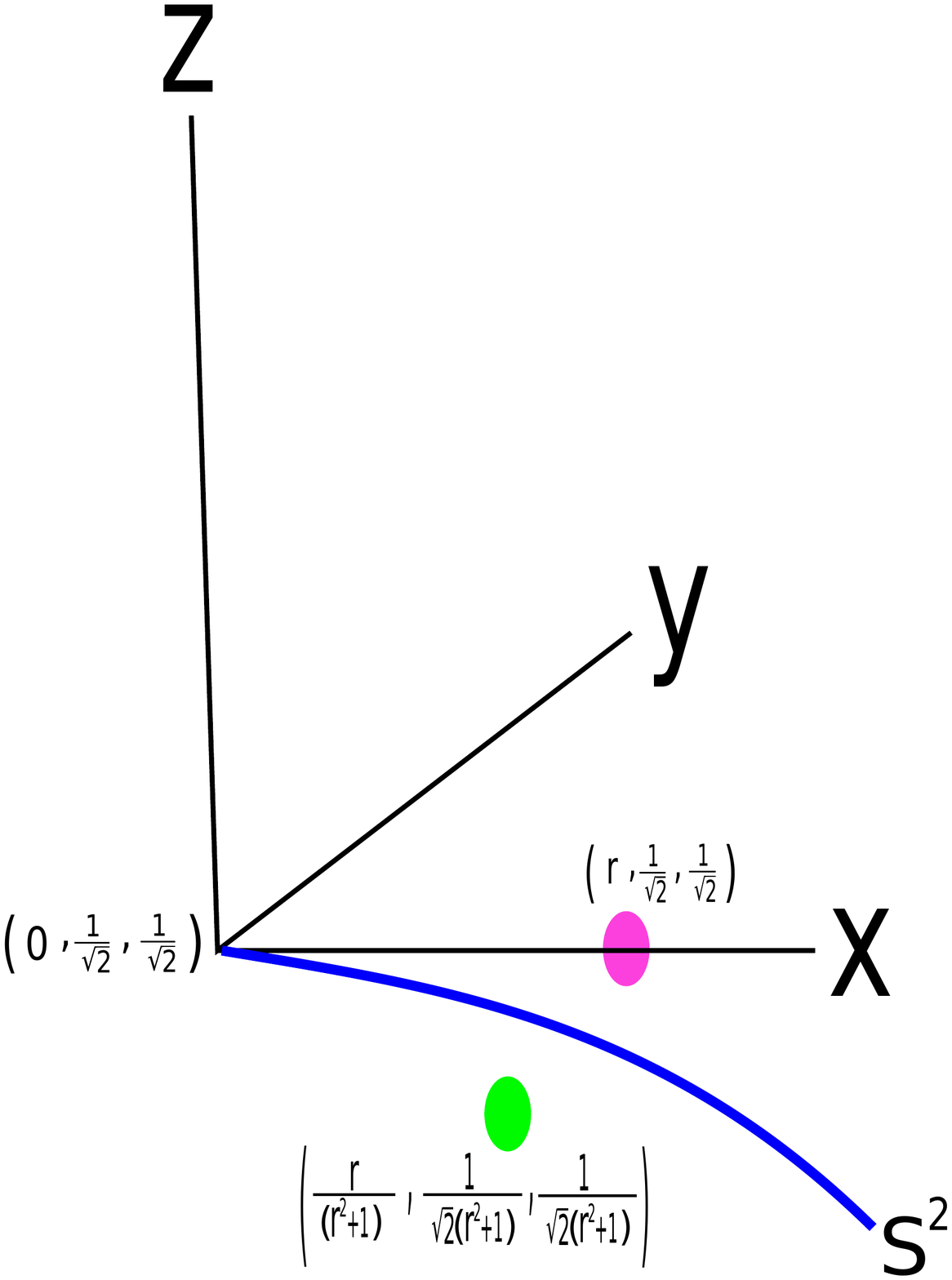}}

The green point represents the image, and the blue curve represents a piece of $S^2$ that is intersecting the attaching sphere $E$. If we now apply
the antipodal map to the above green point $(\frac{r}{r^2 + 1},\frac{1}{\sqrt{2}(r^2 + 1)}, \frac{1}{\sqrt{2}(r^2 + 1)})$ we can see that
we will get a point in $E'$ that will look like what we had before (see two diagrams above).

We now want to carry out the same analysis for a piece of a parallel curve that hits the attaching sphere $E$. We start by choosing our
component of parallel curve to lie directly below the pink 2-handle component. The following diagram shows a close up of this

\centerline {\graphicspath{ {framing/} } \includegraphics[width=3cm, height=3cm]{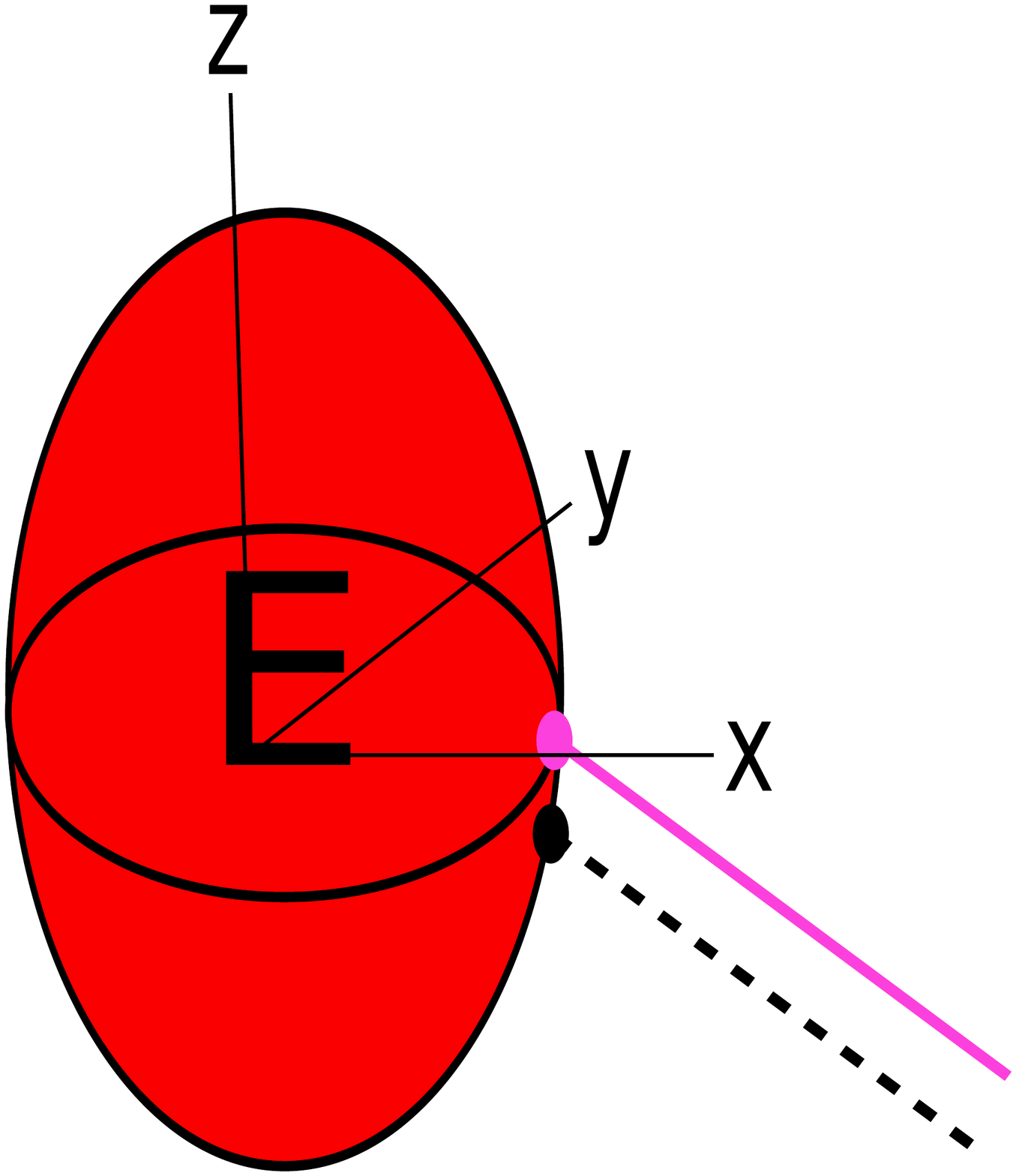}}

We choose the parallel curve so that it hits the attaching sphere at a point with co-ordinates of the form
$(r - \gamma, \frac{1}{\sqrt{2}}, \frac{1}{\sqrt{2}} - \delta)$, where  $\gamma$ and $\delta$ are small so that $(r - \gamma, \frac{1}{\sqrt{2}}, \frac{1}{\sqrt{2}} - 
\delta)$  lies outside of $S^2$ (i.e. has norm greater than one). When we apply inversion in $S^2$ the point
$(r - \gamma, \frac{1}{\sqrt{2}}, \frac{1}{\sqrt{2}} - \delta)$ will map to 
$\bigg(\frac{r-\gamma}{(r-\gamma)^2 + \frac{1}{2} + (\frac{1}{\sqrt{2}}-\delta)^2}, 
 \frac{1}{\sqrt{2}} \cdot \frac{r-\gamma}{(r-\gamma)^2 + \frac{1}{2} + (\frac{1}{\sqrt{2}}-\delta)^2}, 
 (\frac{1}{\sqrt{2}} - \delta) \cdot \frac{r-\gamma}{(r-\gamma)^2 + \frac{1}{2} + (\frac{1}{\sqrt{2}}-\delta)^2}\bigg) $.
The following diagram shows the images of these points under inversion in $S^2$.

\centerline {\graphicspath{ {framing/} } \includegraphics[width=6cm, height=4.5cm]{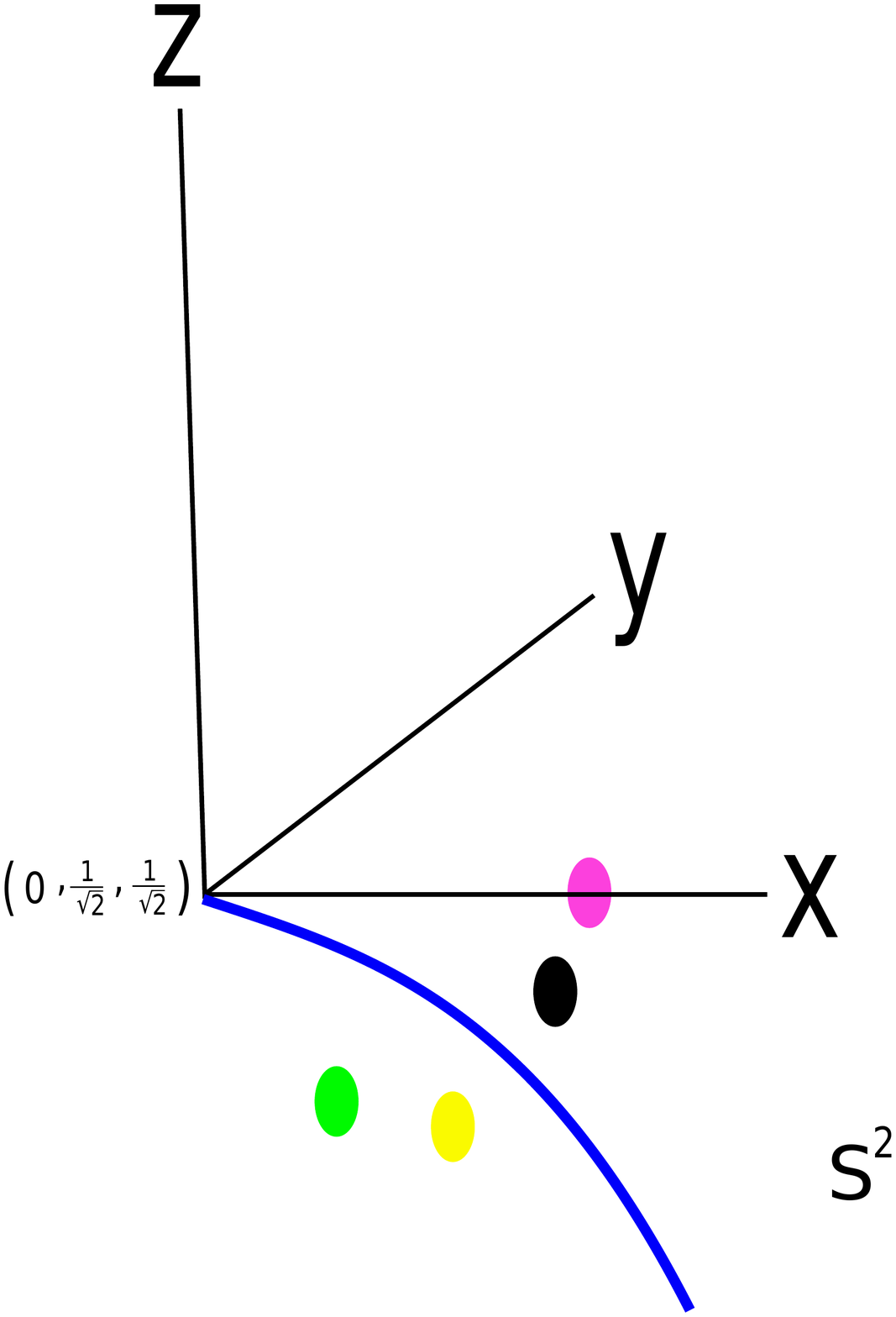}}

The pink and green points are the same points as before, the black point is $(r - \gamma, \frac{1}{\sqrt{2}}, \frac{1}{\sqrt{2}} - \delta)$, and
the yellow point represents the image of the black point under the inversion map.
We then apply the antipodal map to find the complete image of the black point.
The following diagram shows this image, together with how the parallel curve comes out of $E'$. The point to take away is that the parallel curve comes
out of $E'$ lying a bit below the pink 2-handle component.

\centerline {\graphicspath{ {framing/} } \includegraphics[width=3.5cm, height=3.5cm]{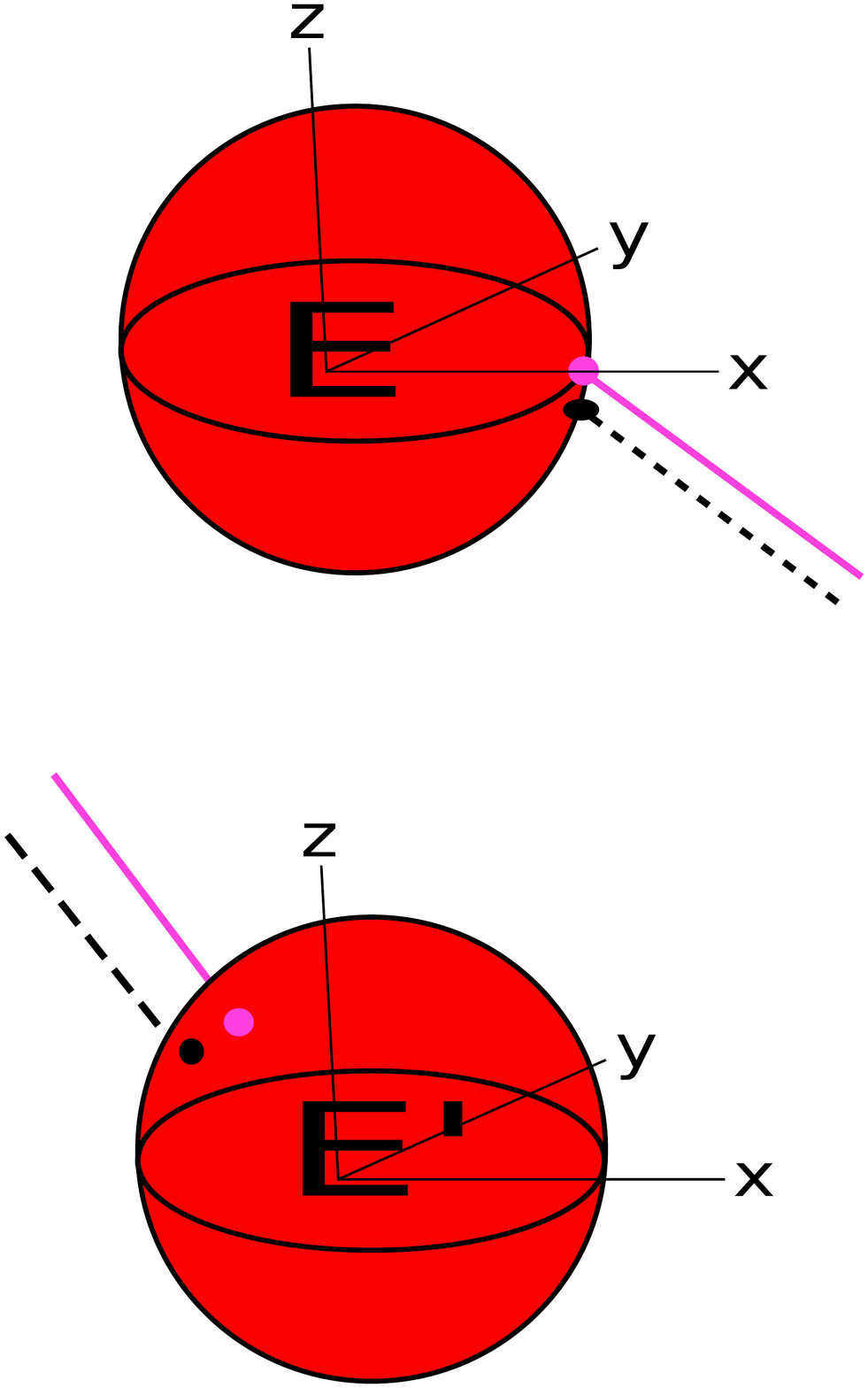}}

The pink 2-handle component then runs and hits the attaching sphere $D'$, along with the parallel curve this looks like.

\centerline {\graphicspath{ {framing/} } \includegraphics[width=6cm, height=6cm]{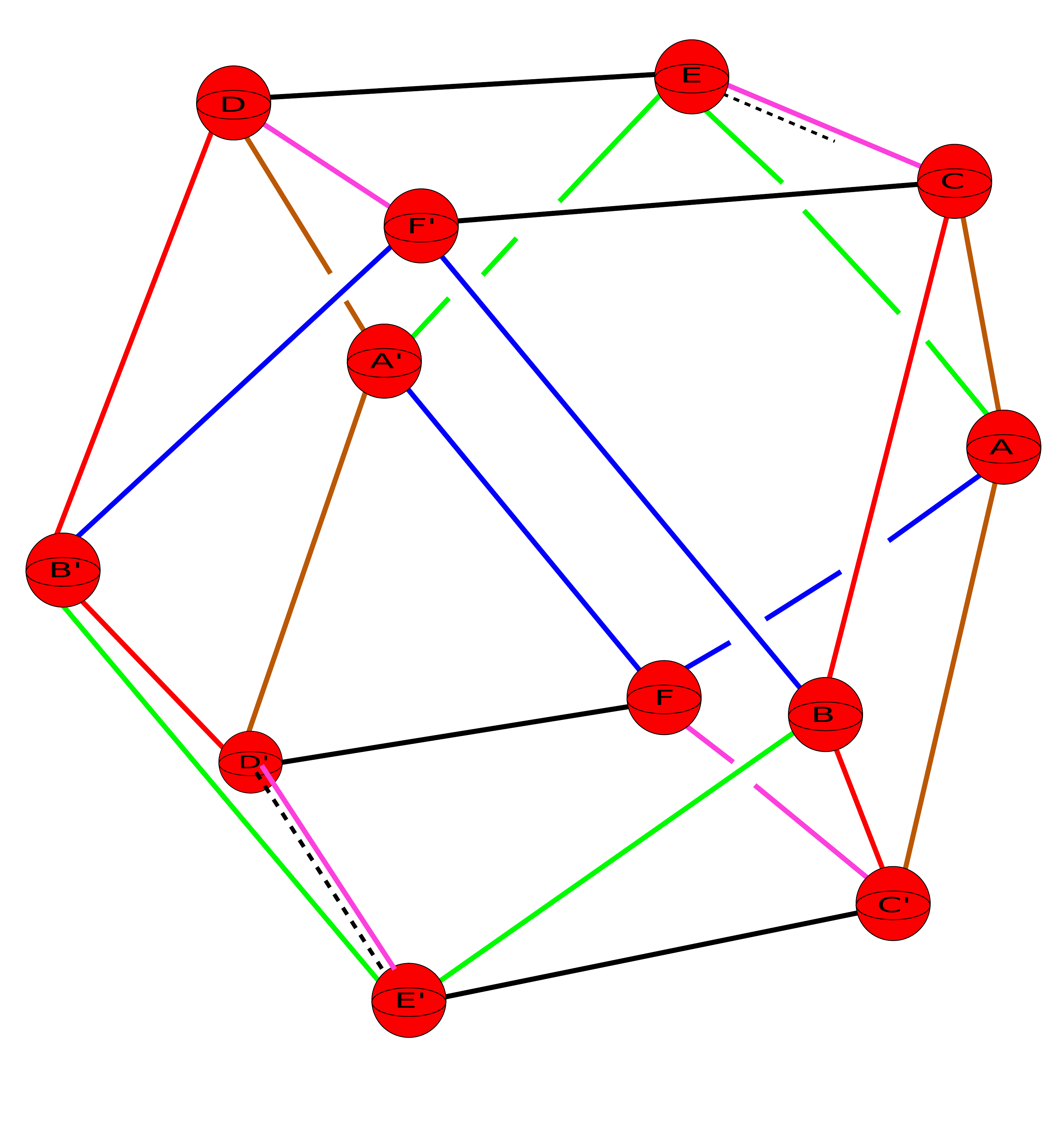}}

The attaching map for the 1-handle $D,D'$ is given by the reflection $(x,y,z) \mapsto (x,y,-z)$, and in this case it easy to see how the
parallel curve component comes out of $D$. The parallel curve will run slightly above the pink 2-handle component.

\centerline {\graphicspath{ {framing/} } \includegraphics[width=6cm, height=6cm]{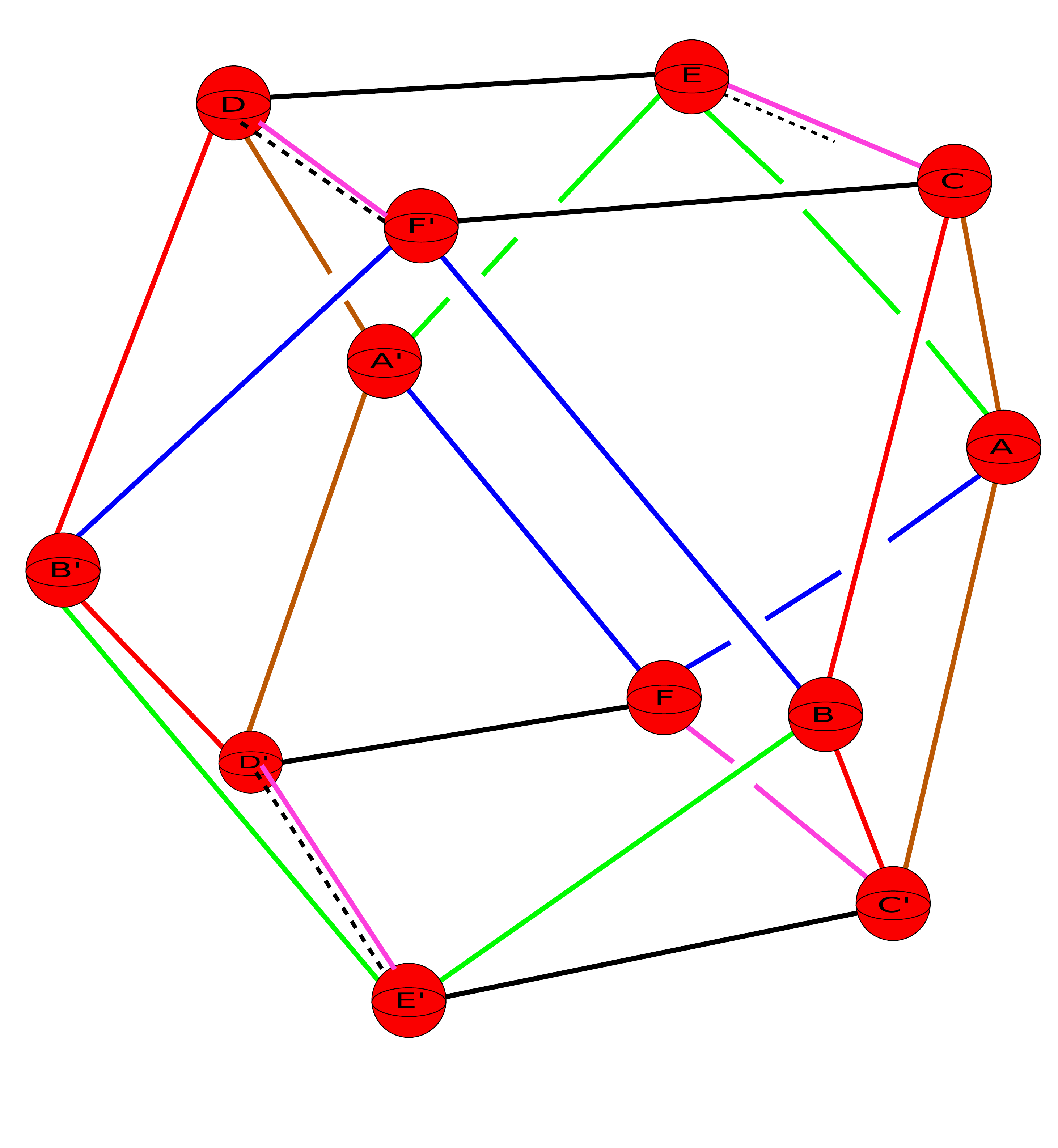}}

The attaching map for $F, F'$ is the same as the attaching map for $E, E'$, and in order to see how the parallel curve component comes out of
$F$ we need to perform a similar sort of analysis as we did for when we tried to understand how it came out of $E'$.  As the details are exactly
analogous to what we did for $E,E'$ we won't give the details. The parallel curve comes out a little above the pink 2-handle
component as shown in the following diagram. \\

\centerline {\graphicspath{ {framing/} } \includegraphics[width=6cm, height=6cm]{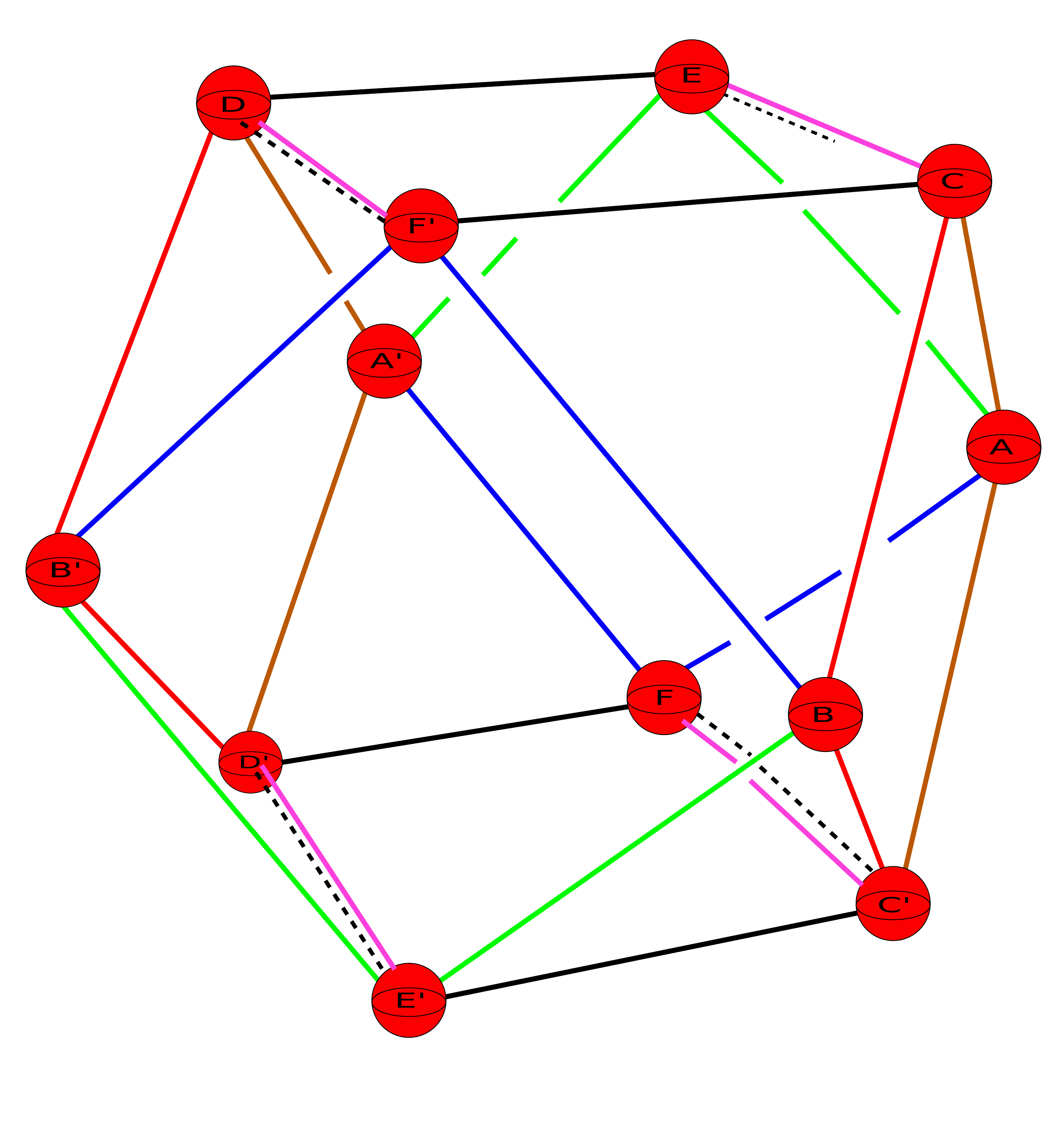}}

Finally, the attaching map for $C, C'$ is the same as that of $D, D'$, namely the reflection $(x,y,z) \mapsto (x,y,-z)$. It is easy
to see that under this map the parallel curve component comes out of $C$ just below the pink 2-handle component and joins up with where our
parallel component initially started going into $E$. The following picture shows this.

\centerline {\graphicspath{ {framing/} } \includegraphics[width=6cm, height=6cm]{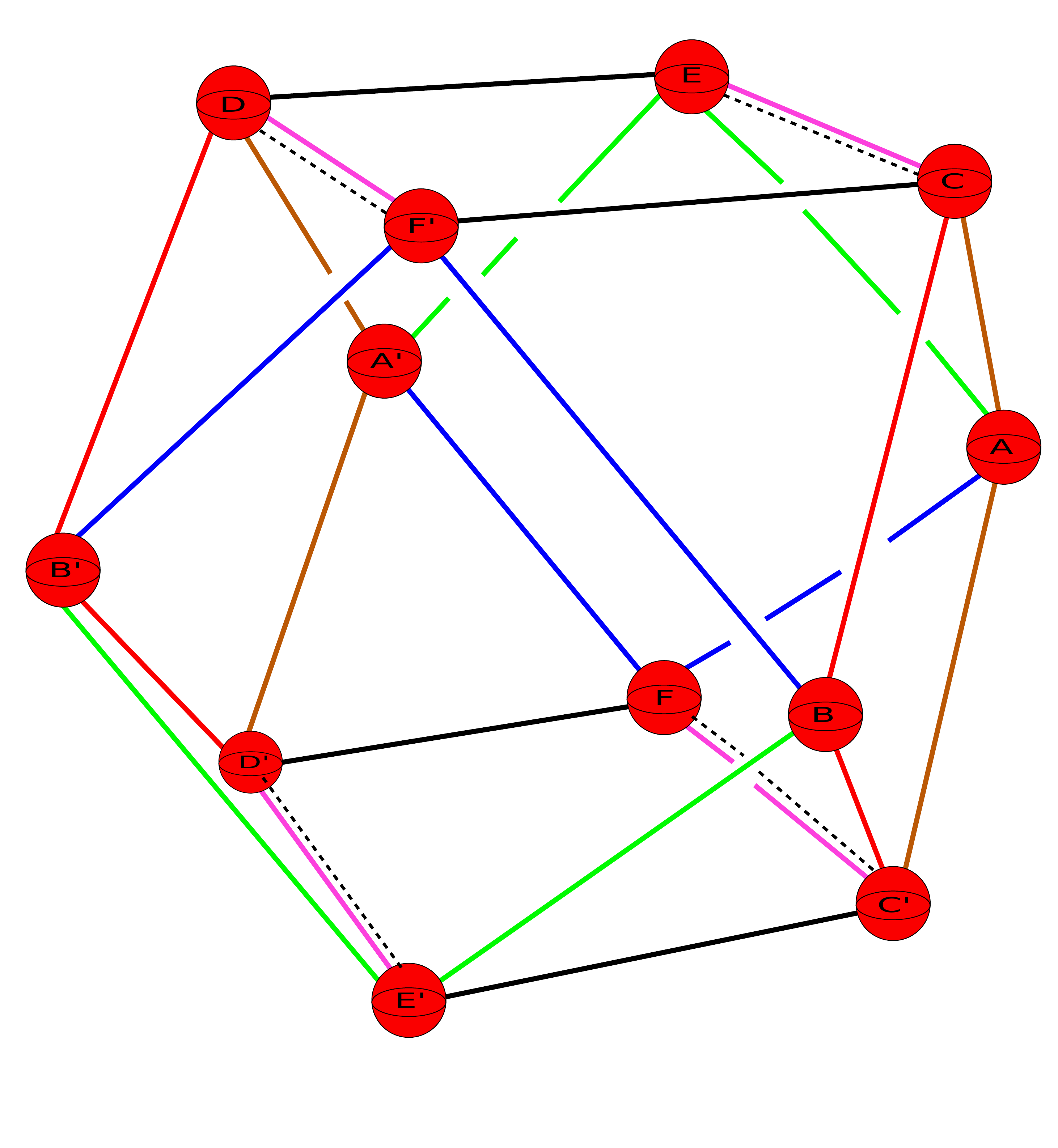}}

What we have seen is that when we take a parallel curve to the pink 2-handle component and follow it around, nowhere does it cross over the pink 2-handle. 
In some sense a parallel curve behaves as if the pink 2-handle was confined to a 2-plane. Another way to think of this is that the normal
bundle can be given by taking a normal vector to each straight line segment comprising a 2-handle and then using the Euclidean connection, given
by the Euclidean dot product, to parallel transport the normal vector along the 2-handle segment. This viewpoint makes it clear that these 2-handles
have parallel curves that behave very similar to the 2-handles that lie in the $x-y$, $x-z$ or $y-z$ planes. Due to this we say that
each 2-handle contained in the Kirby diagram of manifold 1011 has a \textit{planar framing}. The word \textit{planar} is used to remind the reader
that the parallel curves are behaving as if they were lying in a single 2-plane. 


Although we have focused on manifold 1011 a similar sort of analysis can be applied to show that any one of the Ratcliffe-Tschantz manifolds has a Kirby
diagram where all the 2-handles have this \textit{planar framing} type behaviour.

\end{document}